\documentstyle[12pt]{article}

\def\mathrm{{}}
\def\textrm{{}}
\def\a{\alpha}
\def\b{\beta}
\def\C{{\bf C}}
\def\PM{{\cal{PM}}}
\def\H{{\cal{H}}}

\def\c{\gamma}
\def\d{\delta}

\def\M{{\cal{M}}}

\def\N{{\cal{N}}}
\def\P{{P}}
\def\pf{{\noindent\bf Proof: }} 

\def\s{\sigma}
\def\S{\Sigma}
\def\t{\tau}
\def\Z{{\bf Z}}
\def\ZZ{{\cal{Z}}}
\def\BR{{\bf R}}
\def\square{{\vcenter{\hrule height.4pt
      \hbox{\vrule width.4pt height5pt \hskip5pt
	   \vrule width.4pt}
      \hrule height.4pt}}}
\def\qed{\hfill$\square$}
\def\semidirect{\mathop
 {\rlap{\kern 7pt\vrule height 5pt width .4pt depth 0pt} \times}}
\def\Ker{{\rm Ker}}
\def\Im{{\rm Im}}
\newtheorem{proposition}{Proposition}[section]
\newtheorem{theorem}[proposition]{Theorem}
\newtheorem{corollary}[proposition]{Corollary}

\title{Geometric subgroups of mapping class groups}
\author{Luis Paris \and Dale Rolfsen}

\begin{document}

\maketitle

\abstract{This paper is a study of the subgroups of mapping class groups of
Riemann surfaces, called ``geometric'' subgroups, corresponding to the
inclusion of  subsurfaces.  Our analysis includes surfaces with boundary and with
punctures. The centres of all the mapping class groups are calculated. We
determine the kernel of inclusion-induced maps of the mapping class group of a
subsurface, and give necessary and sufficient conditions for injectivity. In the
injective case, 
we show that the commensurability class of a geometric
subgroup completely determines up to isotopy the defining subsurface, and
we characterize centralizers, normalizers, and 
commensurators
of geometric subgroups. }

\bigskip
\centerline{{\it Mathematics Subject Classification:} 
Primary 57N05; Secondary 20F38}

\section{Introduction}

Throughout the paper, $M$ will denote a  {\it compact, connected, oriented}
surface.   The boundary $\partial M$, if nonempty, is a finite collection of
simple closed curves. Consider a finite subset $\P = \{p_1, \dots, p_m\}$, of
$m$ distinct points (often called ``punctures'' or ``marked points'') in the 
interior of $M$.  Define $\H(M,\P)$ to be the group of orientation-preserving
homeomorphisms $h:M \to M$ such that $h$ is the identity on each boundary
component of $M$ and $h(\P) = \P$.  Our main object of study is the {\it mapping
class group} $\M(M,\P) = \pi_0(\H(M,\P)),$ the set of isotopy classes of these
mappings, with composition as the group operation.  We emphasize that throughout
an isotopy, the boundary, and also the points $\P$ remain fixed.  It is clear
that, up to isomorphism, these groups do not depend on the choice of $\P$, but
depend only on the cardinality $m = |\P|$, so we may write $(M,m)$ or $\M(M,m)$ 
in place of $(M,P)$ or $\M(M,P)$,  and simply $\M(M)$ for $\M(M,\emptyset)$. 
$\M(M,P)$ may equivalently  be considered as the group of orientation-preserving
{\it diffeomorphisms} of
$(M,P)$, up to smooth isotopy. 
We refer the reader to the survey articles \cite{Bi2}, 
\cite{Iv1}, \cite{Iv2},
\cite{Mi} and \cite{Mo} and their bibliographies for more information.

Let $N \subset M$ be a {\it subsurface}, by which we mean a closed subset which
is also a surface and for which we always assume the further properties: (1)
every component of $\partial N$ 
lies in the interior of $M$, (2) $P \cap \partial N = \emptyset$.

The inclusion $i:(N,N\cap \P)  \to  (M, \P)$ induces a natural mapping
$$i_* : \M(N,N\cap \P) \to \M(M,\P).$$
If $[h]$ is a class of a homeomorphism of $N$, then $i_*([h])$ is represented
by extending $h$ to $M$ using the identity mapping on $M \setminus N$.
The image $i_*(\M(N,N\cap \P))$ will be called a {\it geometric subgroup} of
$\M(M,\P).$  

Our study of these subgroups depends on a careful analysis of
curves in $M$ and Dehn twists, which are the subject of  Section 3. 
The mapping $i_*$ is often, but not always, injective.   
We determine its kernel in Section 4. Section 5 is devoted
to the centres of mapping class groups. These are
certainly well-known to specialists and many of them
can be found in the literature (see \cite{Iv2} and \cite{IM}).
However, we need the general result for the remainder of the
paper and the proofs are straightforward applications of
Section 3. Our main result is that, assuming injectivity
of the $i_*$, up to a finite number of exceptions, two
geometric subgroups are commensurable if and only if they
are equal if and only if their respective defining
subsurfaces are isotopic (Theorem 6.5). Note that the
assumption that $\partial N$ lies entirely in the interior
of $M$ is necessary for the conclusion of Theorem 6.5;
indeed, without this assumption, it is very easy to
construct non-isotopic subsurfaces $N$ and $N'$ which
define the same geometric subgroup, $N$ satisfying
$\partial N\cap \partial M=\emptyset$ and $N'$ satisfying
$\partial N'\cap\partial M\neq\emptyset$. From the main
result, still assuming injectivity of the $i_*$, we
characterize the commensurator, the normalizer and the
centralizer of a geometric subgroup in $\M(M,P)$.

We close this introduction and
illustrate the injectivity question by  discussing some basic examples, which
are well-known (c.f. \cite{Bi1} and \cite{FLP}).

\bigskip
\noindent{\bf Examples:} (1) $\M(D^2) \cong \{ 1\}$ and $\M(D^2,1) \cong \{
1\}$, where
$D^2$ is a disk.

(2) Similarly, for the 2-sphere, $\M(S^2)$ and $\M(S^2,1)$ are trivial.

The above examples are essentially the only surface mapping class groups which
are trivial.

(3) $\M(S^1 \times I) \cong \Z$.

The mapping class group of the annulus $S^1 \times I$ is generated by a {\it Dehn
twist}, described in Section 3.

(4) $\M(S^1 \times I, 1) \cong \Z^2$.  Here, Dehn twists along the two 
boundary components are generators.  They are not isotopic because of the
puncture.

(5)  As a family of examples, consider disks $D_1 \subset D_2 \subset \dots$ in
the complex plane, where the diameter of $D_m$ is taken to be the real line
interval
\break
$[1-m/(m+1), m + 1/2].$  
Take $P_m = D_m \cap \Z =\{1,\dots, m\}.$ Then it is also
well-known that the mapping class groups $\M(D_m,P_m)$ are isomorphic with the
classical braid groups $B_m$ of Artin, which have generators $\s_1,\dots,
\s_{m-1}$ subject to the relations

$$\s_i\s_j = \s_j\s_i,~~ |i-j| > 1, \hskip2cm \s_i\s_{i+1}\s_i =
\s_{i+1}\s_i\s_{i+1}.$$
Under the isomorphism, $\s_j,~~ 0 < j < m$ corresponds with the (class of the)
diffeomorphism consisting of a `` half-twist'' interchanging the integers $j$
and $j+1$, and supported on a small neighborhood of the interval $[j, j+1]
\subset \C.$  See \cite{Bi1} and \cite{FLP} 
for details, but beware some
differences in choice of conventions.  It is classical, but nontrivial, that for
$n < m$ the homomorphism $B_n \to B_m$, taking $\s_j \in B_n$ to
$\s_j \in B_m$ is injective, allowing us to write $B_n \subset B_m$. We chose
notation so that under the isomorphisms, $B_n \to B_m$ corresponds to the
inclusion-induced mapping $i_* :\M(D_n,P_n) \to \M(D_m, P_m)$. We conclude that
in this case $i_*$ is injective.  Note that the closure of the complementary
subsurface is an annulus with $m-n$ punctures.
The commensurator, the normalizer and the centralizer of
$\M(D_n,P_n)$ in $\M(D_m,P_m)$ are characterized in
\cite{FRZ} and \cite{Ro}.

(6) The following is an example of an inclusion map which is not injective on
mapping class groups.   Take $M = S^2$ a 2-sphere with $P =$ 2 points in $S^2$,
and let
$D$ be a disk in $S^2$ which encloses the points $P$.  We have the map
$i_* :\M(D,P) \to \M(S^2, P)$.  As already discussed, $\M(D,P)$ is the braid
group $B_2$, which is infinite cyclic, generated by
$\s_1$.  However, $\s_1^2$ is isotopic with a Dehn twist along
$\partial D$.  In the larger surface $M = S^2$ this twist is isotopic with the
identity (rel $P$).   So the kernel of $i_*$ in this case is the infinite cyclic
subgroup of index 2 in $\M(D,P)$. $\M(S^2, 2)$ is cyclic of order 2.

(7) For the torus $T^2 = S^1 \times S^1$ with either zero or one puncture,  the
mapping class group is the modular group of invertible $2 \times 2$ matrices
with integer entries:
$\M(T^2) \cong \M(T^2,1) \cong SL(2,\Z)$.  Dehn twists along the curves
$S^1 \times *$ and $*\times S^1$ generate the mapping class group.  Note that 
if $A$ is an annulus neighborhood of one of these curves, and happens to 
enclose the puncture of $\M(T^2,1)$, the map $i_*: \M(A,1) \to \M(T^2,1)$  
fails to be injective.   

\section{Subgroups of mapping class groups}

In this section we review some of the literature regarding subgroups of
mapping class groups.  Although there is little published on the geometric 
subgroups, which are the main concern of the present paper, certain other
subgroups are quite well understood.  First, we recall some general
properties of the mapping class groups themselves.  
A more complete survey can be found in ~\cite{Mi}.  
   
For the closed surface $M_g$ of genus $g$, the mapping class 
groups $\M(M_g,m)$
are known to be finitely presented
~\cite{Ge}, \cite{HT},~\cite{Li}
~\cite{Lo}, ~\cite{Mat}, ~\cite{McCo},
~\cite{Wa1}, ~\cite{Wa2}. 
The generators can be taken to be Dehn twists (discussed below)
along curves and half-twists along arcs connecting the punctures.

According to Grossman~\cite{Gr} and Ivanov~\cite{Iv2}, 
$\M(M_g,m)$ is residually finite -- for
every nontrivial element, there is a homomorphism of the mapping class group
onto a finite group which does not kill that element.  
This implies that $\M(M_g,m)$ is Hopfian~(\cite{LS},
Chapter IV) -- every
epimorphism $\M(M_g,m)\to\M(M_g,m)$ is an isomorphism.
Conversely, in a recent paper~\cite{IM} Ivanov and McCarthy
proved that $\M(M_g,m)$ is co-Hopfian -- every monomorphism
$\M(M_g,m)\to\M(M_g,m)$ is an isomorphism.

Although, it contains
torsion elements, $\M(M_g)$ has a finite index subgroup which is
torsion free (see~\cite{Mi} for a sketch of a proof.)
It has recently been shown by Mosher~\cite{Mos} that the mapping class
groups are {\it automatic}.  This implies that the word problem is solvable
(in quadratic time) and many other consequences~\cite{ECHLPT}.

The outer automorphism group of $\M(M_g,m)$ has been determined
by Ivanov~\cite{Iv3} and McCarthy~\cite{Mac2}. It is equal to
$\Z/2\Z$ under the assumptions $g\neq 0$, $m\ge 3$ if $g=1$,
and $m\ge 1$ if $g=2$. 

The abelianization of $\M(M_g)$ is cyclic of order 12, when $g=1$, cyclic of
order 10 when $g=2$, and trivial (that is, $\M(M_g)$ is perfect) for $g>2$
~\cite{Po}.  

Finite subgroups of $\M(M_g)$ have been extensively studied.  The so-called
Nielsen realization problem~\cite{Zi} was solved by
Kerckhoff~\cite{Ke}.  It asserts that for any finite subgroup $F$ of
$\M(M_g)$, there is a complex structure on $M_g$ such that $F$ is realized as
a group of holomorphic automorphisms of $M_g$.  According to a classical
result of Hurwitz~\cite{Hu}, the orders of finite subgroups are bounded:
$|F| \le 84(g-1), ~g > 1$.

McCarthy~\cite{McCa} showed that subgroups of $\M(M_g)$ satisfy the Tits
alternative: every subgroup either contains a free group on two generators, or
a solvable subgroup of finite index.  Birman, Lubotzky and McCarthy~\cite{BLM}
proved that solvable subgroups of $\M(M_g)$ are virtually abelian, and gave
upper bounds for the rank of free abelian subgroups. 
Ivanov~\cite{Iv2} proved these two results for $\M(M_g,m)$ 
and showed that: if $G$ is a subgroup of $\M(M_g,m)$ which is
not virtually abelian, then $G$ contains
an uncountably infinite number of 
maximal subgroups of infinite index. 
The {\it Frattini subgroup} $\phi(G)$ of a group $G$
is the intersection of all its maximal subgroups.
Ivanov also proved that the
Frattini subgroup of a finitely generated subgroup of
$\M(M_g,m)$ is nilpotent.
Note that this property also holds for
finitely generated linear groups \cite{Plato}.

The centres of the $\M(M_g,m)$ are well-known~\cite{Iv2},
~\cite{IM}: cyclic of order two if $g=1$ and $m\le 2$, and if
$g=2$ and $m=0$, and trivial otherwise.
We determine the centres in the more
general situation, with boundary, in the present paper. 

The Torelli subgroup of the mapping class group consists of classes
of mappings which induce the identity on the homology of the surface. 
These subgroups have been studied extensively in the series of papers
\cite{Johnson1}, \cite{Johnson2}, \cite{Johnson3}, \cite{Johnson4}.
In particular, they describe a finite set of generators for the Torelli
groups. 
The cohomological properties of mapping class groups 
and the Torelli groups
have been studied
intensively and are also well-described in ~\cite{Mi} and
\cite{Mo}.  We shall not use
any of these deeper properties -- indeed our methods are quite elementary and
self-contained, requiring as background only some basic properties of curves on
surfaces due to Epstein ~\cite{Ep}. 
 
\section{Curves and Dehn twists}

Working within the context of a surface $M$ with punctures $P$ as described
above, we shall consider a {\it simple closed curve} in $M \setminus P$  as an
embedding $c: S^1 \to M \setminus P$
which does not intersect the boundary of $M$.   
Note that $c$ has an orientation; the
curve with the opposite orientation, but same image will be denoted $c^{-1}$. 
By abuse of notation, we also use the symbol $c$ to denote the image of $c$.  
We will say that $c$ is {\it essential} if it does not  bound a disk in $M$
disjoint from $\P$, and that $c$ is {\it generic} if  it does not bound a disk
in $M$ containing 0 or 1 point of $\P$.

Two simple closed curves $a, b$ are {\it isotopic} if there exists a continuous 
family $h_t\in\H(M,P),~~ t \in [0, 1]$ of homeomorphisms
such that $h_0$ is the identity and $h_1 \circ a = b$.
Isotopy of curves is an equivalence relation
which we denote by $a \simeq b$.  Following
\cite{FLP} the  {\it index of intersection} of two simple closed curves $a$ and $b$ is:

$$I(a,b) = \inf\{|a' \cap b'|; a' \simeq a, b' \simeq b \}$$

We note that:

1) $I(a,b) = \inf\{|a' \cap b|; a' \simeq a \}$;

2) If $a$ is not generic, then $I(a,b) = 0$ for every simple closed curve $b$;

3) If $a \simeq b$ then $I(a,b)=0$.

A {\it bigon} cobounded by two simple closed curves $a$ and $b$ in 
$M\setminus P$ is a disk 
$D \subset M\setminus P$ 
whose boundary is the union of an arc of $a$ and an arc of $b$.

\begin{proposition}[{Epstein \cite{Ep}}] \label{p1} Let $a, b: S^1 \to M \setminus P$
be two essential simple closed curves, and suppose $a$ is isotopic to $b$.

i) If $a \cap b = \emptyset$, then there exists an annulus in $M \setminus P$
whose boundary components are $a$ and $b$.

ii) If $a \cap b \ne \emptyset$, and they intersect transversely, then
$a$ and $b$ cobound a bigon. \qed

\end{proposition}

\begin{proposition} \label{p2}

Let $a, b: S^1 \to M \setminus P$ be two essential simple closed curves, which 
intersect transversely.  Then
$$I(a,b) = |a \cap b|$$ if and only if $a$ and $b$ do not cobound a bigon.
\end{proposition}

\pf  It is clear that if $a$ and $b$ cobound a bigon, one can isotop 
one of the curves across
the bigon and reduce the cardinality of the intersection 
by two.  Now suppose they do not
cobound a bigon, and choose a simple closed curve $a'$ isotopic to $a$  and
transverse to $a$ such that
$$|a' \cap b| = I(a,b).$$ We will argue by induction on $|a' \cap a|$ that
$$|a \cap b| = |a' \cap b| = I(a,b).$$

\begin{figure}[ht]
\centerline{\input{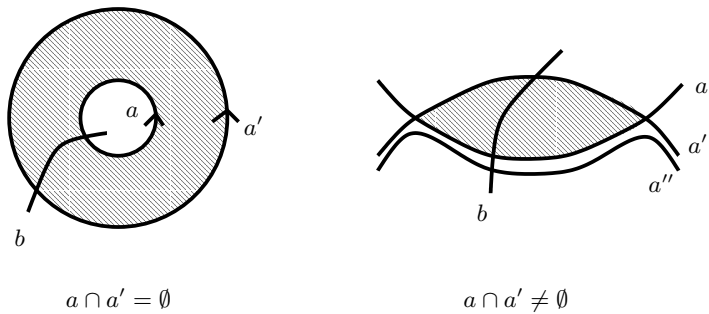}}
\caption{\label{f1} Curves cobounding an annulus and a bigon.}
\end{figure}

If $|a' \cap a| = 0$, then by Proposition \ref{p1}, there is an annulus in $M
\setminus P$ with boundary components $a$ and $a'$.  Each arc of intersection of
$b$ with the  annulus must run from one boundary component to the other
(see Figure \ref{f1}), since
neither $a,b$ nor $a',b$ cobound bigons.
Therefore $|a \cap b| = |a' \cap b| =
I(a,b).$

Now suppose $|a' \cap a| > 0$.  By Proposition \ref{p1}, $a$ and $a'$ cobound a bigon. We
may assume $a' \cap a \cap b$ is empty, so any arc of intersection of $b$  with
the bigon must have one endpoint in $a$ and the other in $a'$
(see Figure \ref{f1}), again because
neither $a,b$ nor $a',b$ cobound bigons.
Therefore we may push $a'$ across the
bigon to obtain
a new curve $a''$ isotopic to $a$ and satisfying:

$$|a'' \cap a| = |a' \cap a| - 2 \quad {\rm and} \quad
 |a'' \cap b| = |a' \cap b| = I(a,b).$$ By inductive hypothesis
$|a \cap b| = |a'' \cap b| = I(a,b).$ \qed

\noindent
{\bf Definition:}  If we parametrize $S^1$ as the unit circle in the
complex plane, and the interval $I = [0,1]$, then the prototype {\it Dehn twist} 
$\t : S^1 \times I \to S^1 \times I $ is
given by

$$\t(z,t) := (ze^{2\pi it},t).$$

Note that $\t$ is the identity on the boundary circles.  More generally, let 
$a: S^1 \to M \setminus P$ be a simple closed curve, and let $N \subset M
\setminus P$ be an annulus regular neighborhood of the image of $a$,
parametrized by 
$\tilde{a}: S^1 \times I \to N$.  Define the {\it Dehn twist along $a$} to be
(the isotopy class of) the homeomorphism
$A(x) =  \tilde{a} \tau \tilde{a}^{-1}(x)$ for $x \in N, ~~A(x) = x$ for $x$
outside $N$.

\begin{figure}[ht]
\centerline{\input{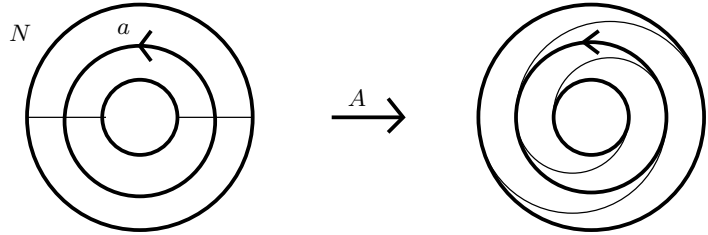}}
\caption{\label{f3} Dehn twist along curve $a$.}
\end{figure}

We will use the convention throughout that a curve is denoted by a lower case
letter and a Dehn twist along the curve is denoted by the corresponding
capital letter.  Note that, depending on parametrization chosen, there are two
choices of Dehn twist along $a$, inverse to each other.  Usually, the choice is
immaterial, provided one is consistent throughout, but we will adopt the
convention in illustrations that a curve crossing $a$ will make a ``right turn"
at each encounter with $a$, after being acted on by $A$ (see Figure \ref{f3}).
We also observe:

1) The Dehn twist along $a^{-1}$ coincides with the Dehn twist along $a$.

2) The curve $a$ is fixed by the Dehn twist $A$.

3) If two curves are isotopic, then so are their corresponding Dehn twists.

4) If $h$ is a homeomorphism of $M$, the Dehn twist along $h(a)$ is $hAh^{-1}$.

5) If $a$ is not generic, $A$ is isotopic to the identity.

\begin{proposition} \label{p3} Let $a,b: S^1 \to M \setminus P$ be two simple
closed curves, $A$ the  Dehn twist along $a$ and $n$ any integer.  Then 
$$I(A^n(b),b) = |n|\cdot I(a,b)^2.$$

\end{proposition}

\begin{figure}[ht]
\centerline{\input{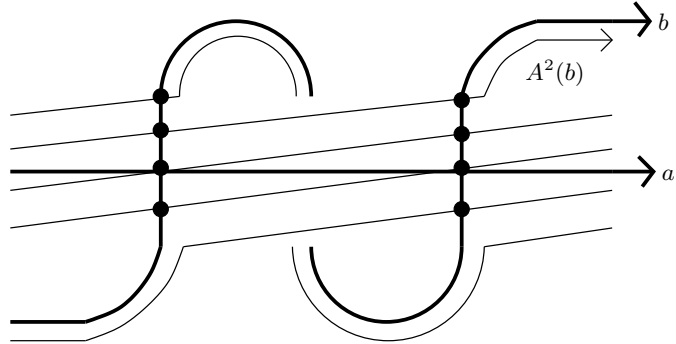}}
\caption{\label{f4} Intersection of $b$ with $A^2(b)$.}
\end{figure}

\pf  This is a special case of a formula in \cite{FLP}.  We outline a  proof, leaving
details to the reader.  Assume $|a \cap b| = I(a,b)$.   The cases $n = 0$ or
$I(a,b) = 0$ being trivial, suppose  they are nonzero.  Construct the curve
$A^n(b)$, which can be seen to cross $b$ exactly $|n|I(a,b)$ times at each point
of  intersection of $a$ with $b$ (see Figure \ref{f4}).
The proof is completed by noting that this is
the minimal intersection of $A^n(b)$ with $b$, up to isotopy.  For otherwise, by
Proposition \ref{p2} there would be a bigon cobounded by them.   One can see
this would imply that $a$ and $b$ also cobound a bigon, which is impossible,
again by Proposition \ref{p2}. \qed

\begin{proposition} \label{p4}  Suppose $a_1, \dots , a_p: S^1 \to M \setminus
P$   are generic simple closed curves such that:

a) $a_i \cap a_j = \emptyset$ if $i \ne j$,

b) $a_i$ is neither isotopic with $a_j$ nor $a_j^{-1}$, if $i \ne j$,

c) none of the $a_i$ is isotopic with a boundary component of $M$.

Then for each $i$, $1 \le i \le p$, there exists a simple closed curve 
$b: S^1  \to M \setminus P$ such that $a_j \cap b = \emptyset$ if $i \ne j$, and
$|a_i \cap b| = I(a_i,b) > 0$.

\end{proposition}

\noindent
{\bf Remark:} The last condition implies that $b$ must be generic in $M \setminus P$.

\bigskip
\pf  We ``cut open'' $M$ along all the curves $a_i$ to obtain the  connected
compact surfaces $N_1, \dots, N_r$ with the property that the union of the
interiors of the $N_j$ is the interior of 
$M \setminus \bigcup_{i=1}^p a_i$ and each boundary component of $N_j$ is either
a boundary component of $M$ or  a copy of some $a_i$.  There is a continuous
projection of the disjoint union onto $M$:
$$\rho : \coprod_{j=1}^r N_j \to M,$$ which covers each curve $a_i$ twice, and
is injective on the union of the  interiors of the $N_j$.  Now fix $i \in
\{1,\dots,r\}$ and consider the curve 
$a_i = \rho(c_1) = \rho(c_2)$, where $c_1$ is a component of the  boundary of
some $N_j$ and $c_2$ is a component of the boundary of $N_k$.

\begin{figure}[ht]
\centerline{\input{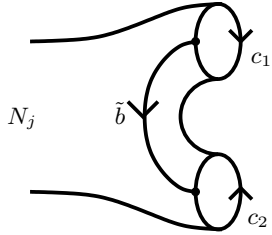}}
\caption{\label{f5} Constructing $b$, case 1.}
\end{figure}

Case 1, $j=k$:  Then $c_1$ and $c_2$ are different components 
of the boundary of the connected
surface $N_j$.  There is an arc $\tilde{b}$ in $N_j$ with one endpoint in $c_1$, 
the other endpoint being the point of $c_2$ having same projection in $a_i$, 
and the interior of $\tilde{b}$ in the interior of $N_j$ and avoiding  the
puncture set $P$ (see Figure \ref{f5}).
Then take $b = \rho(\tilde{b})$.  Clearly 
$a_j \cap b = \emptyset$ if $i \ne j$; moreover    
$a_i$ and $b$ are transverse and they cannot cobound a bigon, so 
$|a_i \cap b| = 1 = I(a_i,b)$.

\begin{figure}[ht]
\centerline{\input{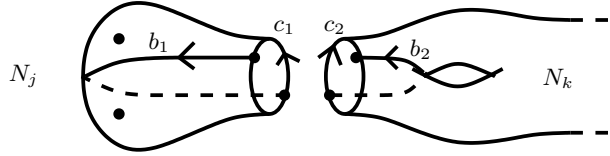}}
\caption{\label{f6} Constructing $b$, case 2.}
\end{figure}

\bigskip
Case 2, $j \ne k$:  Then $c_1$ is a boundary component of $N_j$, and by 
hypothesis, if $N_j$ is a disk, then $N_j \cap P$ contains at least two points, 
and if $N_j$ is an annulus, $N_j \cap P$ contains at least one point.  So in any 
case, there exists an arc $b_1$ in $N_j$ with both endpoints in $c_1$, interior
in the interior of $N_j$ and disjoint from $P$, and such that $b_1$ and $c_1$ do
not cobound a bigon in $N_j \setminus P \cap N_j$ (see Figure \ref{f6}).
In the same way we  choose
an arc $b_2$ in $N_k$, whose endpoints in $c_2$ project to 
$\rho(b_1 \cap c_1)$ in $a_i$, whose interior is interior to 
$N_k \setminus P \cap N_k$, and such that $b_2$ and $c_2$ do not cobound a bigon
in $N_k \setminus P \cap N_k$.  Define $b = \rho(b_1 \cup b_2)$.  Then
$a_j \cap b = \emptyset$ if $i \ne j$ and $|a_i \cap b| = I(a_i,b) = 2$, 
because we have arranged that $a_i$ and $b$ do not cobound a bigon in 
$M \setminus P$. \qed 

\bigskip
Now consider a subsurface $N \subset M$; recall this includes the assumption 
that $P \cap \partial N = \emptyset$ and $\partial N$ is
interior to $M$.  We will say that $N$  is
{\it essential} if 
each component of $\overline{M \setminus N}$ which is
a disk has nonempty  intersection with the puncture set P.

\begin{figure}[ht]
\centerline{\input{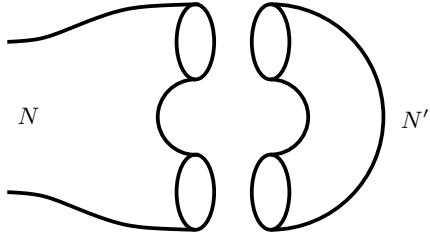}}
\caption{\label{f7} $N'$ is an exterior cylinder for $N$.}
\end{figure}

A component $N'$ of $\overline{M \setminus N}$ 
will be called an {\it exterior cylinder} if
$N'$ is a cylinder (= annulus) disjoint from $P$, with both components of 
$\partial N'$ also being components of $\partial N$
(see Figure \ref{f7}).

\begin{proposition} \label{p5}

Let $N \subset M$ be an essential subsurface and let 
$a, b: S^1 \to N \setminus N \cap P$ be essential simple closed curves. Assume
that $a$ is not isotopic in $N \setminus N \cap P$ to a boundary  component of
an exterior cylinder.  Then $a$ and $b$ are isotopic in 
$M \setminus P$ if and only if they are isotopic in $N \setminus N \cap P$.

\end{proposition}

\pf  The only nontrivial part is to show that if $a$ and $b$ are  isotopic in $M
\setminus P$, then they are isotopic in $N \setminus N \cap P$. We assume $a$
and $b$ intersect transversely and argue by induction on
$|a \cap b|$.

If $|a \cap b| = 0$, then Proposition \ref{p1} implies there exists an annulus
in $M \setminus P$ with (unoriented) boundary $a \cup b$.  Since $a$ is not
isotopic to a 
boundary component of an exterior cylinder
and $N$ is essential, the annulus is
disjoint from $\partial N$, and therefore it lies in $N \setminus N \cap P$. It
follows that  $a$ and $b$ are  isotopic in $N \setminus N \cap P$.

If  $|a \cap b| > 0$, then $a$ and $b$ cobound a bigon in $M \setminus P$, by
Proposition \ref{p1}. Because $N$ is essential, the bigon is disjoint from
$\partial N$, and  therefore it lies in $N$.  Pushing across this bigon defines
an isotopy in
$N \setminus P$ from $b$ to a curve $b'$ with $|b' \cap a| = |b \cap a| - 2$. By
inductive hypothesis, $b'$ is isotopic with $a$ in $N \setminus N \cap P$, so
the same is true of $b$. \qed

\begin{proposition} \label{p6}

Consider two generic simple closed curves $a, b: S^1 \to M \setminus P$,  and
let $A$ and $B$, respectively, denote Dehn twists along these curves. If $j$ and
$k$ are integers, $j \ne 0$, such that $A^j = B^k$ in $\M(M,P)$,
then $a$ is isotopic to $b$ or $b^{-1}$ in $M \setminus P$.

\end{proposition}

\noindent
{\bf Remark:} This proposition and the next one can be found in
\cite{IM} for $\partial M=\emptyset$. The general case needs
an extra argument to consider Dehn twists 
along boundary curves.

\bigskip
\pf  We will argue that if $a$ is not isotopic with $b^{\pm 1}$  then $A^j \ne
B^k$.  It may be assumed that $a$ and $b$ meet  transversely and that $|a \cap
b| = I(a,b)$.  First assume $I(a,b) > 0$.  Then, by Proposition \ref{p3},
$$I(A^j(b), b) = |j|I(a,b)^2 > 0$$
$$I(B^k(b), b) = I(b,b) = 0$$ and we conclude $A^j \ne B^k$.

Now suppose that $I(a,b) = 0$.  If $M$ has nonempty boundary, consider the 
larger closed surface $\hat{M}$ obtained by gluing a torus minus a disk to each
boundary component of $M$ (if $\partial M = \emptyset$, let 
$\hat{M} = M$). By
Proposition \ref{p5} $a$ is not isotopic with $b^{\pm 1}$ in 
$\hat{M}\setminus P$.  By Proposition \ref{p4} there is a simple closed curve
$c$ in $\hat{M}\setminus P$ such that $b \cap c = \emptyset$ and 
$|a \cap c| = I(a,c) > 0$.  Then 
$$I(A^j(c),c) = |j|I(a,c)^2 > 0$$
$$I(B^k(c),c) = I(c,c) = 0,$$ and therefore $A^j \ne B^k$ in $\M(\hat{M},P)$; so
$A^j \ne B^k$ in $\M(M,P)$. \qed

\begin{proposition} \label{p7}

Consider two generic simple closed curves $a, b: S^1 \to M \setminus P$,  and
let $A$ and $B$, respectively, denote Dehn twists along these curves. If $j$ and
$k$ are integers, $j \ne 0$ and $k\neq 0$, 
such that $A^j$ and $B^k$ commute in $\M(M,P)$,
then $I(a,b) = 0$.

\end{proposition}

\pf  Assuming $A^j$ and $B^k$ commute, we have
$$A^j = B^k A^j B^{-k} = C^j,$$ where $C$ is the Dehn twist along the curve $c =
B^k(a)$.  By  Proposition \ref{p6} it follows that $c$ is isotopic with $a^{\pm
1}$. Proposition \ref{p3} implies
$$0 = I(a^{\pm 1},a) = I(c,a) = I(B^k(a),a) = |k|I(a,b)^2,$$ so $I(a,b) = 0$. \qed

\bigskip\noindent
{\bf Remark:} 
An alternative proof of Proposition 3.7 can be
deduced from \cite{Is}, where it is shown that
$A$ and $B$ commute
if $I(a,b)=0$, $A$ and $B$ satisfy the braid relation
$ABA=BAB$ if $I(a,b)=1$, and $A$ and $B$ generate a free
group if $I(a,b)\ge 2$.

\begin{proposition} \label{p8}

Suppose $a_1, \dots , a_p: S^1 \to M \setminus P$   are generic simple closed
curves which are pairwise disjoint, and no curve $a_i$ is isotopic to $a_j$ or
$a_j^{-1}$, $i \ne j$.  Consider the function $$h:\Z ^p \to \M(M,P)$$ defined by
$$h(n_1, \dots,n_p) = A_1^{n_1} \cdots A_p^{n_p},$$ where $A_i$ is the Dehn
twist about $a_i$. Then $h$ is an injective homomorphism.

\end{proposition}

\pf Because the curves are disjoint, the Dehn twists commute, and
$h$ is a homomorphism.  To see it is injective, suppose
$A_1^{n_1} \cdots A_p^{n_p}$ is the identity of $\M(M,P)$ for some
$(n_1, \dots,n_p)$.  We again employ the trick of considering the closed surface
$\hat{M}$, which is $M$ plus a copy of a torus minus a disk glued to each
boundary component.  Clearly, each $a_i: S^1 \to \hat{M}$ is generic, and by
Proposition
\ref{p5} $a_i$ is not isotopic to $a_j^{\pm 1}$ in $\hat{M} \setminus P$,  when
$i \ne j$. Now fix an index $i \in \{1,\dots,p\}$.  Proposition \ref{p4}
supplies a simple  closed curve $b$ in $\hat{M}\setminus P$ disjoint from $a_j,
~~ i \ne j$, with
$$|a_i \cap b| = I(a_i,b) > 0.$$ We calculate, using commutativity of the
twists, 
$b \cap a_j = \emptyset, ~~i \ne j$, and Proposition \ref{p3}:
$$0 = I(b,b) = I(A_1^{n_1} \cdots A_p^{n_p}(b),b)$$
$$ = I(A_i^{n_i}(b),b) = |n_i|I(a_i,b)^2.$$ Therefore $n_i = 0$ and we have
shown $i_*\circ h$ is injective, where 
$i:(M,P) \to (\hat{M},P)$, so $h$ is injective. \qed

\begin{corollary} \label{c9}

If $a: S^1 \to M \setminus P$ is a generic simple closed curve, then the  Dehn
twist about $a$ has infinite order in  $\M(M,P)$. \qed

\end{corollary}

\begin{proposition} \label{p10}

Let $a_1, \dots , a_p, b_1, \dots , b_p: S^1 \to M \setminus P$ be essential
simple closed curves satisfying:

1) $a_i \cap a_j = \emptyset$ and  $b_i \cap b_j = \emptyset$ if $i \ne j$;

2) $a_i$ is not isotopic with $a_j^{\pm 1}$ and 
$b_i$ is not isotopic with $b_j^{\pm 1}$ if $i \ne j$;

3) $a_i$ is isotopic to $b_i$ for each $i = 1,\dots,p$.

Then there exists an isotopy $h_t\in\H(M,P)$ such that
$h_0 = id$ and $h_1\circ a_i = b_i$ for all $i = 1,\dots,p$.

\end{proposition}

\pf We will use a double induction.  First, induction on $p$.  The proposition is
obvious if $p = 1$, so we assume it is true for $p-1$ pairs of curves.   This
means that, replacing each $a_i$ by $h_1\circ a_i$, we may assume that $a_i =
b_i$ for 
$i = 1, \dots, p-1$.  Then we have $a_p$ disjoint from $a_j = b_j, ~~j<p$ and 
also $b_p$ disjoint from $a_j = b_j, ~~j<p$, and $a_p$ isotopic in 
$M \setminus P$ to $b_p$.  We will be done if we show that there is a further
isotopy taking $a_p$ to $b_p$ which does not move the curves $a_j = b_j,
~~j<p$.  Taking
$a_p$ and $b_p$ to be transverse, we will argue by induction on $|a_p \cap b_p|$.

If $|a_p \cap b_p| = 0$, then $a_p$ and $b_p$ cobound an  annulus in $M
\setminus P$, by Proposition \ref{p1}. Any simple closed curve in this annulus
must be either inessential or parallel to  a boundary component, so our
hypotheses guarantee that the annulus is disjoint  from all the curves $a_j =
b_j, ~~j<p$.  Then there is an isotopy across the  annulus taking $a_p$ to
$b_p$; the isotopy may be taken to be the identity  outside a small neighborhood
of the annulus, so the other curves do not move. Suppose $|a_p \cap b_p| > 0$,
then by Proposition \ref{p1}, the curves cobound  a bigon in $M \setminus P$,
and we may argue as above that the bigon is  disjoint from the other curves.  An
isotopy taking $a_p$ across the bigon,  fixed outside a neighborhood of the
bigon, reduces the number $|a_p \cap b_p|$ and does not move the other curves. 
The inductive hypothesis now gives a  final isotopy taking $a_p$ to $b_p$. \qed

\section{Subsurfaces and injectivity}  

Define the {\it pure} mapping class group of $M$ relative to $P$ to be the
subgroup 
$\PM(M,P)$ of $\M(M,P)$ consisting of all (classes of) diffeomorphisms which
fix $P$ pointwise. 
Letting $\S_P$ denote the group of permutations of the set $P$, we
have the exact sequence
$$1 \to \PM(M,P) \to \M(M,P) \to \S_P \to 1.$$

\begin{figure}[ht]
\centerline{\input{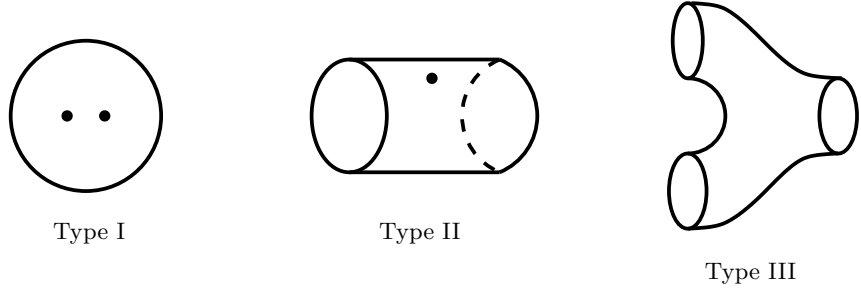}}
\caption{\label{f8} The three types of pantalons.}
\end{figure}

\bigskip\noindent
{\bf Definitions:} A pair equivalent to $(D^2,2)$ is called a {\it pantalon of
type I} (see Figure \ref{f8}).  
We have already observed that $\M(D^2,2)$ is infinite cyclic,
generated by a ``half-twist'' $\s$ which interchanges the two punctures. 
$\PM(D^2,2)$ is the subgroup generated by $\s^2$, which is (up to isotopy) the
same as a Dehn twist along the boundary.
  
A pair $(S^1 \times I, 1)$ is a {\it pantalon of type II} (see Figure \ref{f8}).   
The mapping class
groups  
$\M(S^1 \times I, 1)$ and $\PM(S^1 \times I, 1)$ coincide, and are isomorphic
with $\Z^2$, generated  by the Dehn twists along the two boundary components.
  
A {\it pantalon of type III} is a connected planar surface $M$ with  three
boundary components, the puncture set $P$ is taken to be empty
(see Figure \ref{f8}).  
Its mapping
class group is $\M(M,0) \cong \Z^3$, generated by Dehn twists  along the three
boundary curves.

\begin{figure}[ht]
\centerline{\input{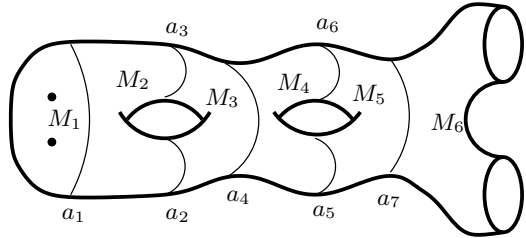}}
\caption{\label{f9} A pantalon decomposition.}
\end{figure}

\bigskip
We now define one of our basic tools, a {\it pantalon decomposition} of a pointed
surface $(M,P)$ (see Figure \ref{f9}).  This consists of:

1) a collection $(M_i,P_i), i = 1, \cdots, r$, each pair being a pantalon of
one of the three types, together with maps $\phi_i: (M_i,P_i) \to (M,P)$,

2) a collection $a_1, \cdots, a_p$ of simple closed curves in $M$, disjoint
from each other, from $P$ and from $\partial M$ satisfying the following.

a) Each $\phi_i$ is injective on the interior of $M_i$.

b) $\phi_i(int M_i)$ and $\phi_j(int M_j)$ are disjoint if $i \ne j$. 

c) $\phi_i$ takes each boundary component of $M_i$ to one of the curves
$a_k$ or to a boundary component of $M$.  Two boundary components of $M_i$ are
allowed to map to the same curve $a_k$.

d) $M = \cup_{i=1}^r \phi_i(M_i)$ and 
$P = \cup_{i=1}^r \phi_i(P_i)$.

\bigskip
Informally, we say that $a_1, \cdots, a_p$ determine a pantalon decomposition of
$(M,P)$, cutting $M$ open into pantalons $(M_i, P_i)$.

   It is straightforward to check that every connected compact orientable
surface $M$ with punctures $P$ admits a pantalon decomposition, with the
following exceptions:

a) $M = S^2$, a sphere, and $|P| \le 3$;

b) $M = D^2$, a disk, $|P| \le 1$;

c) $M = S^1 \times I$, an annulus, with $P = \emptyset$; 

d) $M = T^2$ the torus, with $P = \emptyset$.

Recall that a subsurface $N \subset M$ is essential if no component of 
$\overline{M \setminus N}$ is a disk disjoint from $P$.   If some component $N'$
of $\overline{M \setminus N}$ is a disk with
$|N' \cap P| = 1$, we call $N'$ a {\it pointed disk} exterior to $N$.

\begin{figure}[ht]
\centerline{\input{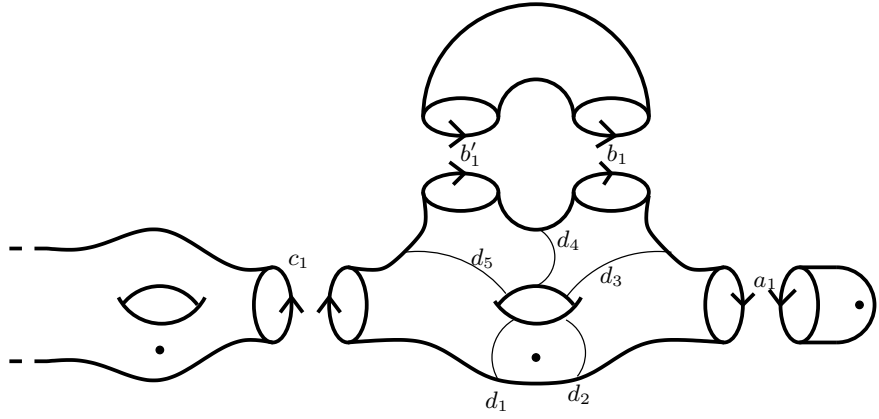}}
\caption{\label{f10} Subsurface $N$ and pantalon decomposition.}
\end{figure}

\begin{theorem} \label{p11}

Consider an essential subsurface $N \subset M$ with the associated  homomorphism
induced by inclusion:
$$i_*: \M(N,N \cap P) \to \M(M,P)$$

i) If $N$ is a disk and $|N \cap P| \le 1$ then $\M(N,N \cap P)$ is trivial, and
therefore $i_*$ is injective.

ii) Suppose $N$ is an annulus and $N \cap P = \emptyset$.  If $N$ has an
exterior pointed disk, then the kernel of $i_*$ is  $\M(N,N \cap P)$; otherwise
$i_*$ is injective.

iii) Assuming that $(N,N \cap P)$ is not as in (i) or (ii), let 
$a_1,\dots,a_r$ denote the boundary components of $N$ which bound pointed disks
exterior to $N$, and let $b_j, b_j', ~ j = 1,\dots,s$ be  the pairs of boundary
components of $N$ which cobound exterior cylinders (disjoint from $P$) 
(see Figure \ref{f10}). Denote
by $A_i, B_j, B_j'$ the Dehn twists corresponding to the curves
$a_i, b_j, b_j'$, respectively.  Then the kernel of $i_*$ is generated by
$$\{ A_1,\dots,A_r , B_1^{-1}B_1', \dots, B_s^{-1}B_s' \}$$ and is isomorphic to
$\Z^{r+s}$.

\end{theorem}

\bigskip
\pf  Part (i) is obvious and (ii) is a direct consequence of Proposition
\ref{p8}.  To prove (iii), let $[h] \in Ker(i_*)$, where $h: N\to N$.  We have
the following commutative diagram, with exact rows.  
$$\matrix{ 1 & \to & \PM(N,P \cap N)& \to & \M(N,P \cap N)  & \to & \S_{P \cap
N} & \to & 1 \cr  &&\downarrow&& i_*\downarrow&&\downarrow&&\cr 1 & \to &
\PM(M,P)& \to & \M(M,P)  & \to & \S_{P} & \to & 1 \cr} $$
 Since the homomorphism $\S_{P \cap N} \to \S_{P}$ is  injective, $[h]$ is in
$\PM(N,P \cap N)$.

Let $c_1, \dots ,c_t$ denote the components of $\partial N$ different from the
$a_i$ and $b_j, b_j'$.  In addition, let 
$d_1, \dots, d_u: S^1 \to N \setminus P$ be simple closed curves which 
determine a pantalon decomposition of $(N,N\cap P)$
(see Figure \ref{f10}).
Note that all the curves
we are considering are pairwise disjoint and non-isotopic.
Since $h$ is isotopic to  the
identity in $M \setminus P$, each $h\circ d_i$ is isotopic to $d_i$ in
$M \setminus P$.  Proposition  \ref{p5} implies that
$h\circ d_i$ is isotopic to $d_i$ in $N \setminus N\cap P$.   
By Proposition \ref{p10}
we may suppose that
$h\circ d_i = d_i$ for all $i = 1,\dots,u$, and that $h$ is the identity on the
boundary of each pantalon.  Using the structure of the  
pure mapping class groups of
pantalons we conclude that 
$$[h] = A_1^{\a_1}\cdots A_r^{\a_r}B_1^{\b_1}B_1'^{\b_1'}\cdots
B_s^{\b_s}B_s'^{\b_s'}C_1^{\c_1}\cdots C_t^{\c_t}D_1^{\d_1}\cdots D_u^{\d_u}.$$
Therefore
$$1 = i_*[h] = B_1^{\b_1 + \b_1'}\cdots B_s^{\b_s + \b_s'} C_1^{\c_1}\cdots
C_t^{\c_t}D_1^{\d_1}\cdots D_u^{\d_u}.$$ By Proposition \ref{p8}
$$\b_1+\b_1' = \cdots = \b_s+\b_s' = \c_1 = \cdots =\c_t = 
\d_1 = \cdots = \d_u = 0.$$ Therefore
$$[h] = A_1^{\a_1}\cdots A_r^{\a_r}(B_1^{-1}B_1')^{\b_1'}\cdots
(B_s^{-1}B_s')^{\b_s'}.$$ Conversely, it is clear that any $[h]$ of this form is
in the kernel of $i_*$. Finally, Proposition \ref{p8} implies that $Ker(i_*)$ 
is isomorphic to $\Z^{r+s}$. \qed

\begin{corollary} \label{p12}

Let $N \subset M$ be any subsurface  (with $\partial N \cap P = \emptyset$) and
let
$$i_*: \M(N,N \cap P) \to \M(M,P)$$ the natural homomorphism.

i) If $N$ is a disk and $|N \cap P| \le 1$, then $i_*$ is injective.

ii) If $N$ is an annulus and $N \cap P = \emptyset$, then $i_*$ is injective if
and only if there is no boundary component of $N$ which is the boundary of a
disk intersecting $P$ in less than two points.

iii) If $(N,N \cap P)$ is not as in (i) or (ii), then $i_*$ is injective  if and
only if no component of $\overline{M \setminus N}$ is 
either an annulus disjoint from $P$ 
whose boundary components are both boundary components of $N$,
or a disk which contains less than two points of $P$.

\end{corollary}

\section{Centres}

\begin{figure}[ht]
\centerline{\input{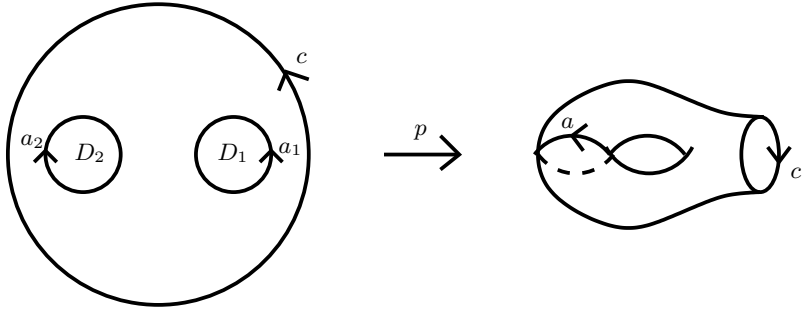}}
\caption{\label{f11} The projection map.}
\end{figure}

We begin this section by considering a special mapping of the surface $M$  of
genus one, and with one boundary component; that is, $M$ is a torus minus  a
disk, with empty puncture set $P$.  We model $M$ as a certain identification
space of a planar surface, as follows (see Figure \ref{f11}).  Let 
$$D = \{ z \in \C ; |z| \le 4 \}$$
$$D_1 = \{ z \in \C ; |z - 2| < 1 \}$$ 
$$D_2 = \{ z \in \C ; |z + 2| < 1 \}$$ 

Then $D \setminus (D_1 \cup D_2)$ is a pantalon (of type III)  with boundary
curves $a_1, a_2, c: S^1 \to \partial M$ which we  parametrize as follows, $0
\le \theta  \le 2\pi$:
$$c(e^{i\theta}) = 4e^{i\theta}$$
$$a_1(e^{i\theta}) = 2 + e^{i\theta}$$
$$a_2(e^{i\theta}) = -2 - e^{-i\theta}.$$

\bigskip
We consider $M = {(D \setminus (D_1 \cup D_2))} / \sim$ where we identify the
points on the curves $a_1$ and $a_2$ by
$$a_1(e^{i\theta}) \sim a_2(e^{i\theta}).$$ Denote the natural projection by 
$$p : D \setminus (D_1 \cup D_2) \to M.$$ The ``meridian'' curve $a : S^1 \to M$
is defined by
$$a = p \circ a_1 = p \circ a_2.$$

Now we define a homeomorphism $R: D \to D$ by the equation:

\begin{displaymath} R(re^{i\theta}) = \left\{ \begin{array}{ll}
   re^{i(\theta - \pi)} & {\rm if}~~ 0 \le r \le 3 \\
   re^{i(\theta - (r-2)\pi)} & {\rm if}~~ 3 \le r \le 4 \\
   \end{array}\right.
\end{displaymath}

We have $R(D_1) = D_2,~~ R(D_2) = D_1, ~~ R\circ a_1 = a_2^{-1}$ and $R\circ a_2
= a_1^{-1}.$
Therefore, $R$ induces a homeomorphism 
$$\tilde{R} : M\to M$$ such that $$\tilde{R}\circ a = a^{-1}.$$ Its class $\rho
= [\tilde{R}] \in \M(M)$ will be called a {\it half-twist} of $M$ 
along $c$ relative to
$a$.  The Dehn twists $A$ and $C$ about the  respective curves $a$ and $c$
relate to $\rho$ as follows:
$$\rho^2 = C, ~~ \rho A \rho^{-1} = A.$$ The latter equation follows since $\rho
A \rho^{-1}$ is the Dehn twist about
$\rho \circ a = a^{-1}$, and the Dehn twist about $a^{-1}$ is the same as the
Dehn twist about $a$ for any simple closed curve $a$.

\begin{proposition} \label{p13}

Let $M, \rho, a$ and $A$ be as in the above discussion.  If $G$ is the subgroup
of 
$\M(M)$ consisting of classes of homomorphisms $h: M \to M$ such that
$h\circ a$ is isotopic with $a$ or $a^{-1}$, then $G$ is generated by 
$\{A ,\rho\}$ and is isomorphic with $\Z^2$.

\end{proposition}

\pf Let $[h] \in G$; we may assume $h\circ a = a$ or $h\circ a = a^{-1}$.   Note
that $a$ determines a pantalon decomposition of $M$, with a single pantalon of
type III, namely $D \setminus (D_1 \cup D_2)$.

If $h\circ a = a$, by the structure of the mapping class group of the pantalon,
for some integers $m,n$:
$$[h] = A^nC^m = A^n\rho^{2m}.$$

If $h\circ a = a^{-1}$, then $\tilde R\circ h\circ a = a$ and $\rho
[h]$ has the form $$\rho [h] = A^nC^m = A^n\rho^{2m}.$$  Therefore
$$[h] = \rho^{-1}A^n\rho^{2m}= A^n\rho^{2m-1}.$$

It is clear that $A$ and $\rho$ are in $G$, so they generate $G$ and we have
already  noted that they commute.  If $A^n\rho^{m} = 1$, then 
$$(A^n\rho^{m})^2 = A^{2n}\rho^{2m} = A^{2n}C^m = 1,$$ and by Proposition
\ref{p8}, $2n = m = 0$.  This shows $G \cong \Z^2$. \qed

\bigskip 
We turn now to determining the centre $\ZZ\M(M,P)$ of an arbitrary 
(compact orientable) surface $M$ with puncture set $P$, that is, the subgroup 
of $\M(M,P)$ consisting of mapping classes which commute with all elements  of
$\M(M,P)$.  First we record some simple cases:
If $M$ is a sphere or disk and $|P| \le 1$, $\ZZ\M(M,P) = 
\M(M,P) = \{1\}$.
$\ZZ\M(S^2,2) = \M(S^2,2) = \Z / 2\Z$.  $\M(S^2,3) \cong \Sigma_3$ and therefore
$\ZZ\M(S^2,3)$ is the trivial group.  Each of the pantalons (type I, II or III) 
has an abelian mapping class group, so the centre is equal to the whole group.
For the same reason $\ZZ\M(S^1 \times I) = \M(S^1 \times I) \cong \Z.$  

\begin{figure}[ht]
\centerline{\input{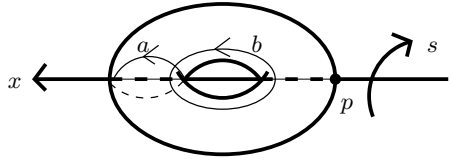}}
\caption{\label{f12} Involution generating the centre of $\M(T^2)$.}
\end{figure}

The case of the torus is somewhat more interesting.  Since 
$$\M(T^2, 0) \cong \M(T^2, 1) \cong SL(2, \Z)$$ we see algebraically that the
centre is the cyclic group of order two, consisting of the two diagonal matrices,
$\pm I$, where $I$ is the identity matrix.  However, to warm up for the more
complicated cases, whe will establish this fact geometrically.   Consider the
torus $T^2$ embedded in $xyz$-space as the set of points of distance 1 from the
circle $x^2 + y^2 = 4,~ z = 0$ (see Figure \ref{f12}).
Let $s: T^2 \to T^2$ be  the
(orientation-preserving) involution $s(x,y,z) = (x,-y,-z)$.  
In the case of
$\M(T^2, 1)$, we suppose that $P = \{p\}$ is a point on the $x$-axis, so that
$s(p) = p$.   Let $\s = [s] \in \M(T^2)$ and $\s_1 = [s] \in \M(T^2, \{p\})$.

\begin{proposition} \label{p14}

The centre of $\M(T^2)$ is the cyclic group of order 2 generated by $\s$; 
similarly $\ZZ\M(T^2, 1) \cong  \Z / 2\Z$, generated by $\s_1$.

\end{proposition}

\pf  We prove the first part; the case of $\ZZ\M(T^2, \{p\})$ being essentially 
the same.  Since $s$ has order 2, $\s^2 = 1$ in $\M(T^2)$.  Let 
$a,b: S^1 \to T^2$ be the circles parametrized by
$$a(e^{i\theta}) = (2 + \cos \theta, 0, \sin \theta),
\quad b(e^{i\theta}) = (\cos
\theta, \sin \theta, 0).$$ Then $$s \circ a = a^{-1} ~~ {\rm and} ~~ s \circ
b = b^{-1}.$$  Since $a$ is not isotopic to $a^{-1}$ we conclude that $\s \ne
1$, and so $\s$ has order 2.  Letting $A, B$ be the Dehn twists about $a, b$ we
have
$$\s A \s^{-1} = A  ~~ {\rm and} ~~ \s B \s^{-1} = B,$$ and since $A$ and $B$
generate  $\M(T^2)$ we conclude that $\s \in \ZZ\M(T^2).$

Now suppose that $[h] \in \ZZ\M(T^2),$ where $h:T^2 \to T^2.$ Then 
$A = [h] A [h]^{-1} = C$ is a Dehn twist about the curve $c = h \circ a.$  By
Proposition \ref{p6}, $c \simeq a^{\pm 1}$, so we may assume that 
$h \circ a = a$ or $h \circ a = a^{-1}$.  Case 1, $h \circ a = a$; we can cut
$T^2$ open along $a$ and conclude from the mapping class of the cylinder that
$[h] = A^k$ for some integer $k$.  Since $I(a,b) = 1$ and $A^k$ and $B$ commute,
Proposition \ref{p7} implies that $k = 0$, and $[h] = 1$ in this case.   Case 2,
$h \circ a = a^{-1}$.  Then $s \circ h \circ a = a$ and we conclude  from case 1
that $\s [h] = 1$, so $[h] = \s^{-1} = \s$. \qed

\bigskip
Next we consider the torus with two marked points.  Using the above model for
$T^2$, let $P = \{(0,3,0), (0,-3,0)\}$, that it, a pair of points which are
interchanged by the involution $s: T^2 \to T^2$ described above.  Let $\s_2 \in
\M(T^2, P) = \M(T^2, 2)$ be the class represented by $s$.  The following is
proved in the same manner as the above.

\begin{proposition} \label{p14a}  The centre of $\M(T^2, 2)$ is the cyclic group
of order 2 generated by $\s_2$.  \qed

\end{proposition}

We now consider the closed surface $M$ of genus 2, which we realize 
in $\BR^3$ as
the boundary of a uniform regular neighborhood of the union of the circles
$$(x+1)^2 + y^2 = 1, ~ z = 0\ {\rm and}\ (x-1)^2 + y^2 = 1, ~ z = 0.$$  
Let 
$s: \BR^3 \to \BR^3$ be the same involution as in the previous discussion, and  let
$\s \in \M(M)$ be represented by the restriction of $s$ to $M$.
 
\begin{proposition} \label{p15}

The centre of $\M(M)$, where $M$ is a closed surface of genus 2, is the
cyclic
group of order 2 generated by $\s$.

\end{proposition}

\begin{figure}[ht]
\centerline{\input{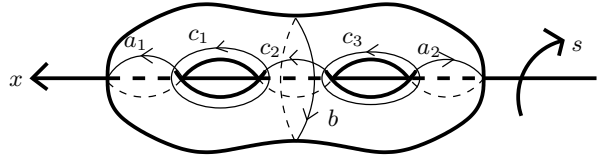}}
\caption{\label{f13} Generator of the centre of $\M(M_2)$.}
\end{figure}

\bigskip
\pf  Let $a_1, a_2, b, c_1, c_2, c_3$ be the simple  closed curves as in
Figure \ref{f13}.
In particular, $b$ is the curve of  intersection of $M$ with the plane
$x=0$,
$c_1$ and $c_3$ lie in the plane $z = 0$  and $a_1, a_2, c_2$ are in the plane
$y=0$. We have $$s \circ a_i = a_i^{-1}, ~~ s \circ b = b, ~~ s \circ c_i =
c_i^{-1}.$$ Since $a_1$ is not isotopic with $a_1^{-1}$, $\s \ne 1$ has order 2
in $\M(M)$. Let $A_1$,$ A_2$,$ B$, $C_1$, $C_2$, $C_3$ be the Dehn twists about
the respective curves, which generate $\M(M)$, according to \cite{Li}.   From the
above equations we conclude that
$$\s A_i \s^{-1} = A_i, ~~ \s B \s^{-1} = B, ~~ \s C_i \s^{-1} = C_i,$$ so $\s$
is indeed central in $\M(M)$.

Given $[h] \in \ZZ\M(M)$ we have $A_i = [h] A_i[h]^{-1}$, which is a Dehn twist 
about $h \circ a_i$.  By Proposition \ref{p6} $h \circ a_i$ is isotopic with
$a_i^{\pm 1}$.  Similarly  $h \circ b$ is isotopic with
$b^{\pm 1}$.  By Proposition \ref{p10} we may assume that $h \circ a_i = a_i$ or
$a_i^{-1}$ and that  $h \circ b = b$ or $b^{-1}.$

Let $M_1$ and $M_2$ be the closures of the two components of $M \setminus b$, 
with $a_1$ in $M_1$, $a_2$ in $M_2$.  Now $h \circ b = b^{-1}$ can only happen
if $h(M_1) = M_2$ and $h(M_2) = M_1$, which is impossible because 
$h \circ a_i = a_i^{\pm 1}$.  So we must have $h \circ b = b$, $h(M_i) = M_i$.
Letting $\rho_i$ denote the half-twist of $M_i$ along $b$ 
relative to $a_i$, Proposition
\ref{p13} implies that $[h]$ can be written in the form
$$[h] = \rho_1^{n_1} A_1^{m_1} \rho_2^{n_2} A_2^{m_2},$$ for some integers
$n_1, n_2, m_1, m_2$. Consider $$ [h]^2 =  A_1^{2m_1} A_2^{2m_2} B^{n_1 + n_2}.$$
Since $[h]^2$ is central it commutes with $C_1$, as do $A_2$ and $B$.  Therefore
$A_1^{2m_1}$ also commutes with $C_1$.  But $I(a_1,c_1) = 1$, so by Proposition
\ref{p7}, $2m_1 = 0$.  Similarly $2m_2 = 0$ and $n_1 + n_2 = 0$.  We note that
either $n_1$ and $n_2$ are both odd or both even.  If they are both even, then
$h \circ a_1 = a_1$ and $h \circ a_2 = a_2$, if both odd, then 
$h \circ a_1 = a_1^{-1}$ and $h \circ a_2 = a_2^{-1}$.

Case 1, $n_1$ and $n_2$ are even: write $n_1 = 2k_1$ and  $n_2 = 2k_2$.  Then
$$[h] = \rho_1^{2k_1} \rho_2^{2k_2} = B^{k_1 + k_2} = B^0 = 1,$$ because $k_1 +
k_2 = (n_1 + n_2)/2 = 0$. 

Case 2, $n_1$ and $n_2$ are odd.  Then
$$s \circ h \circ a_1 = a_1 ~~ {\rm and} ~~ s \circ h \circ a_2 = a_2$$ and
by case 1, $\s [h] = 1$.  Therefore $[h] = \s^{-1} = \s$. \qed

\begin{proposition} \label{p16}

Let $M$ denote the genus one oriented surface with one boundary component.  Then
the centre of $\M(M)$ is the infinite cyclic group 
generated by the
half-twist $\rho$, defined in the discussion preceding Proposition \ref{p13}.

\end{proposition}

\pf Let $a$ be the curve on $M$, $A$ the Dehn twist about $a$, 
$\tilde{R}:M \to M$, $\rho = [\tilde{R}]$, all as described in the discussion
preceding Proposition \ref{p13}.  Let $b$ be the curve on $M$ which is the
image, under the identification $D \setminus (D_1 \cup D_2) \to M$, of the
interval [-1, 1], so that $I(a,b) = 1$.  Let $B$ be the Dehn twist about $b$;
then $\{A, B\}$ generates $\M(M)$.

Noting that $\tilde{R} \circ a = a^{-1}$ and $\tilde{R} \circ b=
b^{-1}$, we see that $\rho$ commutes with $A$ and $B$, and  therefore $\rho$ is
in the centre of $\M(M)$.
Since $\rho^2=C$ and $C$ has infinite order, $\rho$ 
has also infinite order.

Now let $\xi$ be in the centre of $\M(M)$, and $h:M \to M$ a homeomorphism with
$[h] = \xi$.  Then $A = \xi A \xi^{-1}$ is a Dehn twist about $h \circ a$, so by
Proposition \ref{p6} $h \circ a$ is isotopic with $a$ or $a^{-1}$.  By
Proposition \ref{p13} $\xi$ is of the form
$$\xi = A^p \rho^q, \quad p,q \in \Z.$$ Because $\xi$ and $\rho$ commute with
$B$, $A^p$ commutes with $B$.  But $I(a,b) = 1$, so Proposition \ref{p7} implies
that $p = 0$ and therefore $\xi = \rho^q$. \qed

\bigskip\noindent
{\bf Remark:} Let $M$ be a genus one oriented surface with one boundary
component as above. According to \cite{Wa1} $\M(M)$ is isomorphic to the
Artin braid group $\langle A,B|ABA=BAB\rangle$, by \cite{Ch} its center is
the infinite cyclic group generated by $(ABA)^2$, and one can check
directly that $\rho=(ABA)^2$.

\bigskip
Following is the main result regarding the centre of $\M(M,P)$, for a general
Riemann surface.  We have already noted that there are certain exceptional
cases, so we list the hypotheses here for the generic result:

1) If $M$ is a sphere, then assume $|P| \ge 4$,

2) If $M$ is a disk, then assume $|P| \ge 3$,

3) If $M$ is an annulus, then assume $|P| \ge 1$,

4) If $M$ is a torus, then assume $|P| \ge 3$,

5) If $M$ is a surface of genus one and one boundary component, then assume
$|P| \ge 1$,

6) If $M$ is a surface of genus 2, then assume $|P| \ge 1$.

\begin{theorem} \label{p17}

Let $M$ be a connected compact orientable surface with marked points $P \subset
M$ and assume the hypotheses (1) - (6) above.  Let
$c_1, \cdots, c_q: S^1 \to \partial M$ be the boundary curves of $M$ and  
$C_1, \cdots, C_q$ the Dehn twists about these curves.  Then the centre
$\ZZ\M(M,P)$  of $\M(M,P)$ is the subgroup generated by $\{ C_1 \cdots, C_q \}$
and is isomorphic with $\Z^q$.  In particular, if 
$\partial M$ is empty, the centre of $\M(M,P)$ is trivial.

\end{theorem}  

\pf It is clear that $C_i$ is central in $\M(M,P)$; the fact that the subgroup
generated by $C_1, \cdots, C_q$ is isomorphic with $\Z^q$ follows directly from
Proposition \ref{p8}.  Now consider an element $\xi$ in the centre of $\M(M,P)$,
with representative homeomorphism $h:M \to M,~~ [h] = \xi$.

Consider curves $a_1, \cdots, a_p$ which detemine a pantalon decomposition of
$(M,P)$.  (All the cases which do not admit a pantalon decomposition have been
excluded.)
Since $\xi$ is central, the Dehn twist $A_i$ about
the curves $a_i$ satisfy 
$$A_i = \xi A_i \xi^{-1},$$ and since $\xi A_i \xi^{-1}$ is a Dehn twist about
$h \circ a_i$ we conclude from Proposition \ref{p6} that $h \circ a_i$ is
isotopic with $a_i$ or $a_i^{-1}$.  By Proposition \ref{p10} we may suppose that
$$h \circ a_i = a_i^{\pm 1},$$ and consequently, $h$ permutes the pantalons.

\begin{figure}[ht]
\centerline{\input{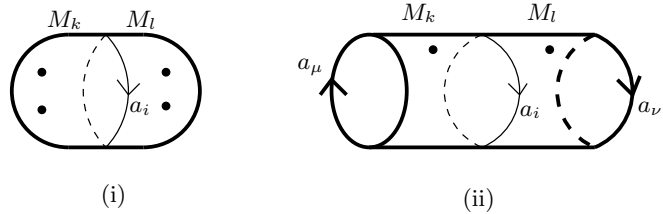}}
\caption{\label{f14} First cases, $a_i$ separating.}
\end{figure}

The pantalons $(M_j,P_j)$ are mapped to $M$ by maps $\phi_j$.  We call a curve
$a_i$ {\it separating} if it is in the image of $\phi_k$ and $\phi_l$, for some
$k \ne l$  
(a separating curve in our sense need not separate 
the surface itself).
First assume that $a_i$ is separating and that 
$h \circ a_i = a_i^{-1}.$  Then $h(\phi_k(M_k)) = \phi_l(M_l)$ and 
$h(\phi_l(M_l)) = \phi_k(M_k)$.  We consider all the possibilities:

If $(M_k,P_k)$ and $(M_l,P_l)$ are pantalons of type I, then $M$ is a 2-sphere
and $|P| = 4$ (see Figure \ref{f14}.i).  Consider the exact sequence:
$$1 \to \PM(S^2, 4) \to \M(S^2, 4) \to \Sigma_4 \to 1.$$ Since the centre of
$\Sigma_4$ is trivial, $\xi$ must act trivially on $P$, whereas $h$ interchanges
the points $\phi_k(P_k)$ and $\phi_l(P_l)$, so this case cannot occur.

If  $(M_k,P_k)$ and $(M_l,P_l)$ are pantalons of type II, then 
$\phi_k(M_k) \cap \phi_l(M_l) = a_i$ and there are other curves $a_\mu$ 
and $a_\nu$
which form the remaining boundary curves of these pantalons (see Figure
\ref{f14}.ii).
Then, since $h$
interchanges the pantalons, $h(a_\mu) = a_\nu^{\pm 1}$.  On the other hand, $h(a_\mu)
= a_\mu^{\pm 1}$, so we must have $a_\mu = a_\nu^{\pm 1}$ and conclude that $M$ is a
torus and $|P| = 2$.  This case has been excluded.

\begin{figure}[ht]
\centerline{\input{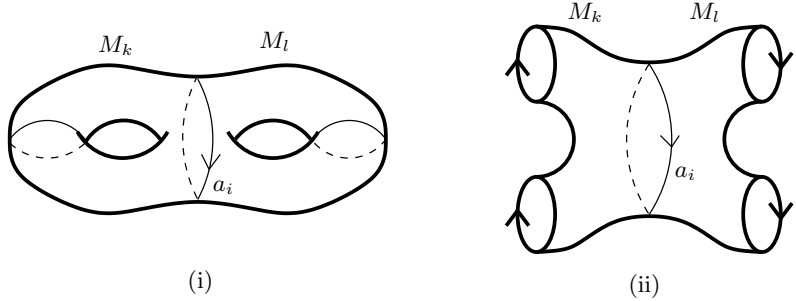}}
\caption{\label{f15} $a_i$ separating type III pantalons.}
\end{figure}

If $(M_k,P_k)$ and $(M_l,P_l)$ are pantalons of type III with identifications,
so that $\phi_k(M_k)$ and $\phi_l(M_l)$ are genus 1 surfaces with one boundary
component $a_i$, then $M$ is a closed surface of genus 2 and $P$ is empty
(see Figure \ref{f15}.i).
This
case also has been excluded.

Finally suppose $(M_k,P_k)$ and $(M_l,P_l)$ are pantalons of type III, mapped
homeomorphically by $\phi_k$ and $\phi_l$ (see Figure \ref{f15}.ii).
Then we argue as in the type II case
that their other boundary components must be identified, and that $M$ is closed,
has genus 2, and $P$ is empty, an excluded case.

\bigskip
Thus we have shown that if $a_i$ is separating, then $h \circ a_i = a_i$.  In
addition we have seen that $h(\phi_k(M_k)) = \phi_k(M_k)$ for all $k$.

\begin{figure}[ht]
\centerline{\input{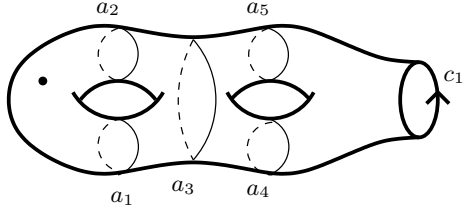}}
\caption{\label{f16} Pantalon decomposition, all $a_i$ separating.}
\end{figure}

Recall the notation $m = |P|$ and $q$ is the number of boundary components of
$M$.  We complete the proof by considering three cases:

Case 1: $m+q \ge 2$.  Then we may assume all the $a_i$ are separating 
(see Figure \ref{f16}). If $M$ is
not a sphere or disk, then each pantalon can be taken to be of type II or III. 
If $M$ is either a sphere or disk, the exact sequence
$$1 \to \PM(M,P) \to \M(M,P) \to \Sigma_P \to 1,$$ and the fact that $\Sigma_P$
has trivial centre under the assumption $|P| \ge 3$, shows that $\xi \in
\PM(M,P)$, or in other words, $h$ fixes $P$ pointwise.  So we have, by the
structure of the (pure) mapping class groups of the pantalons of type I, II and
III:
$$\xi = A_1^{r_1} \cdots A_p^{r_p}C_1^{s_1} \cdots C_q^{s_q},$$ for some
integers $r_1, \dots, r_p, s_1, \dots, s_q.$ Now fix $i \in \{1, \dots, p\}$. 
By Proposition \ref{p4}, there exists a generic simple closed curve $b:S^1 \to
M\setminus P$ such that 
$I(a_i,b) > 0$ but $a_j \cap b = \emptyset$ if $j \ne i$.  If $B$ is the Dehn
twist about $b$, we see that $B$ commutes with $A_j$, $j \ne i$.  $B$ also
commutes with all the $C_j$ and with $\xi$ since they are central.  It follows
that $B$ commutes with $A_i^{r_i}$.  By Proposition \ref{p7} we conclude that
$r_i = 0$.  Since $i$ was arbitrary, $$\xi = C_1^{s_1} \cdots C_q^{s_q}.$$

\begin{figure}[ht]
\centerline{\input{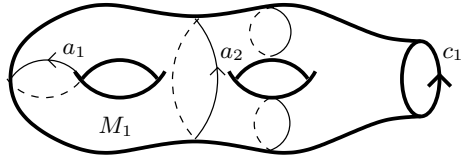}}
\caption{\label{f17} Pantalon decomposition with $a_1$ nonseparating.}
\end{figure}

\bigskip
Case 2: $m+q = 1$.  We may now suppose that $a_1$ is not separating, but $a_2,
\dots, a_p$ are separating curves, and the pantalon $(M_1,P_1)$ is of type III
with two boundary curves identified to $a_1$, so that $\phi_1(M_1)$ is a surface
of genus one with one boundary component, which we may take to be $a_2$
(see Figure \ref{f17}).
Moreover we may assume $(M_2,P_2), \dots ,(M_r,P_r)$ are pantalons of type II or
III, embedded in $M$.  Let $\rho$ denote the half-twist of $\phi_1(M_1)$
along $a_2$
relative to $a_1$.  Then by Proposition \ref{p13} and the structure of the
mapping class groups of type II and III pantalons, $\xi$ can be written
(assuming $q=1, m=0$ and noting that $\rho^2 = A_2$):
$$\xi =\rho^k A_1^{r_1}A_3^{r_3} \cdots A_p^{r_p}C_1^{s_1},$$ for integers
$k,r_1,r_3, \dots, r_p, s_1$.   Then
$$\xi^2 =A_1^{2r_1} A_2^k A_3^{2r_3} \cdots A_p^{2r_p}C_1^{2s_1}.$$  Employing
the same argument as in Case 1, we conclude that 
$$2r_1 = k = 2r_3 = \cdots =2r_p = 0,$$ and so 
$\xi = C_1^{s_1}$ if $q=1$.  If $q=0$ and $m=1$ we similarly conclude that 
$\xi = 1$.

Case 3: $m+q = 0$.  In particular, $q=0$.  Then we need two singular pantalons
in the decomposition of $M$, $\xi$ has an expression as a product of two
half-twists together with powers of the $A_i$ and we conclude by examining
$\xi^2$, exactly as in case 2 that the powers are all zero and therefore $\xi =
1$.   \qed

\section{Commensurability} 

If $G$ is a group, then two subgroups
$H, H' < G$ are said to be {\it commensurable} if $H \cap H'$ has finite index
in both $H$ and $H'$.  Commensurability is an equivalence relation on the set of
all subgroups of $G$, really of interest only for infinite groups.  Following is
an elementary property of commensurable subgroups, which we shall find useful.

\begin{proposition} \label{p18.5}

Suppose $H$ and $H'$ are commensurable subgroups of $G$.  Then for each $h \in
H$, there exists a nonzero integer $k$ such that $h^k \in H'.$

\end{proposition}

\pf If not, then $\{h^k\}, k \in \Z$ is an infinite set of elements, all in
different cosets of $H$ rel $H \cap H'$, contradicting finite index. \qed

\begin{proposition} \label{p19}

Suppose $a_1, \dots, a_p: S^1 \to M \setminus P$ are essential simple closed
curves which are pairwise disjoint.  Let $b: S^1 \to M \setminus P$ be an
essential simple closed curve such that $I(a_i, b) = 0$ for all $i = 1, \dots,
p$.  Then there exists a simple closed curve $c: S^1 \to M \setminus P$ isotopic
to $b$ and such that 
$a_i \cap c = \emptyset$ for all $i=1,\dots,p$. 

\end{proposition}  

\begin{figure}[ht]
\centerline{\input{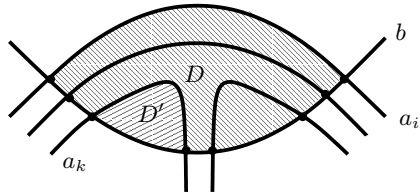}}
\caption{\label{f18} Bigon in proof of Proposition 6.2.}
\end{figure}

\pf  We may assume $b$ transverse to all the $a_i$ and argue by induction on the
cardinality of $b \cap (a_1 \cup \cdots \cup a_p)$.  If 
$b \cap (a_1 \cup \cdots \cup a_p)$ is empty, there is nothing to prove; suppose
it is nonempty.  Choose $i$ so that 
$b \cap a_i \ne \emptyset$.  Proposition \ref{p2} implies that $b$ and $a_i$
cobound a bigon $D$.  There may be other intersections of $a_j$ with $D$, but
there is always an outermost bigon $D' \subset D$ cobounded by $b$ and some
$a_k$, and otherwise disjoint from $a_1 \cup \cdots \cup a_p$
(see Figure \ref{f18}).
Now we can push
$b$ across $D'$ to obtain a curve $b'$ isotopic with $b$ and
$$|b' \cap (a_1 \cup \cdots \cup a_p)| =  |b \cap (a_1 \cup \cdots \cup a_p)| -
2.$$ By inductive hypothesis, there exists an isotopy from $b'$ to an essential
simple closed curve $c$ in $M \setminus P$ such that $a_i \cap c = \emptyset$
for all $i = 1, \dots, p$.  \qed

\bigskip\noindent
{\bf Definitions:}  
Recall that a subsurface is {\it essential} if none of its exterior components
is a disk with zero marked points.
By a {\it marked subsurface} we mean a pair
$(N,Q)$ where $N$ is a subsurface of $M$ and
$Q=P\cap N$.
We will say that an essential marked subsurface $(N,Q)$ 
of $(M,P)$ is {\it
injective} provided it satisfies:

a) If $N$ is a disk, then $|Q| \ge 2$,

b) If $N$ is an annulus with $Q$ empty, then there is no pointed disk component
of the exterior of $N$ in $M$,

c) If $N$ is not an annulus, or if $|Q| \ge 1$, then no component of 
$\overline{M\setminus N}$ 
is a disk with one marked point, or a cylinder with no marked
points and both boundary components in $\partial N$.

By Corollary \ref{p12} these criteria assure the injectivity of  $$\M(N,Q) \to
\M(M,P).$$

If $(N,Q)$ and $(N',Q')$ are marked subsurfaces of  $(M,P)$, 
we say they are {\it
isotopic} provided there is a continuous family of homeomorphisms 
$h_t\in\H(M,P),~~ t\in[0,1]$ such that $h_0 =$ identity,
and $h_1(N,Q) = (N', Q')$.  In particular, $Q = Q'$. 

We would like to be able to say that geometric subgroups are commensurable if
and only if they are equal, if and only if their defining subsurfaces are
isotopic.  However, just as with centres, there are some exceptions to the
general principle.  Our first family are the infinite cyclic geometric
subgroups.  The only marked surfaces with mapping class $\Z$ are 
$(S^1 \times I, 0)$ and $(D^2, 2)$; we note for future reference that these are
the only mapping class groups which contain $\Z$ as a finite-index subgroup.

\begin{figure}[ht]
\centerline{\input{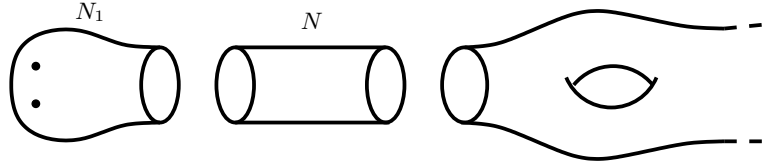}}
\caption{\label{f19} Injective subsurfaces in Proposition 6.3.}
\end{figure}

\begin{proposition} \label{p20}

Suppose $(N,Q)$ and $(N',Q')$ are injective subsurfaces of $(M,P)$, and suppose
$N = S^1 \times I$ and $Q$ is empty. Suppose $\M(N,Q)$ and $\M(N',Q')$ are
commensurable subgroups of 
$\M(M,P)$.  Then either (1) $N'$ is a cylinder $S^1 \times I$ 
with $Q'$ empty or
(2) $N'$ is a disk and $|Q'| = 2$.  

In case (1), $N'$ is isotopic with $N$.  In case (2), 
$(N',Q') = (D^2,2)$, then one of the components $(N_1,Q_1)$ of the exterior of
$(N,Q)$ in $(M,P)$ is a $(D^2,2)$, and 
$(N_1,Q_1)$ is isotopic with $(N',Q')$ (see Figure \ref{f19}).

\end{proposition}  

\bigskip
\pf  Since $\M(N,Q)$ is infinite cyclic, so is every nontrivial subgroup,
including $\M(N,Q) \cap \M(N',Q')$, so $(N',Q')$, 
its mapping class group
containing a finite index $\Z$,
can only be $(S^1\times I, \emptyset)$ or $(D^2,2)$.

Case 1: $(N',Q')$ is $(S^1\times I, \emptyset)$.  Let 
$a'(z) = (z, 1/2),~~ z \in S^1$ denote the central curve of $N'$ and similarly
label the central curve of $N$ as $a: S^1 \to N$, 
and let $A', A$ be Dehn twists of
$M$ about these curves.  $A$ and $A'$ represent generators of $\M(N,Q)$ and
$\M(N',Q')$, respectively.  By commensurability their intersection has finite
index, so there exist nonzero integers $k, l$ such that  $$A^k = A'^{l}.$$  
Proposition \ref{p6} implies that $a'$ is isotopic with $a$ or with $a^{-1}$. 
Since $N$ and $N'$ are regular neighborhoods, respectively, of $a$ and $a'$, it
follows that $N$ and $N'$ are isotopic in this case.

Case 2: $(N',Q')$ is $(D^2,2)$.  Choose notation as in Case 1, except that $a'$
now denotes the boundary curve of $N'$.  Again we have nonzero $k$ and $l$ so
that $A^k = A'^{l}$, and we conclude that $a'$ is isotopic with $a^{\pm 1}$. 
After an isotopy, we may assume $a'$ equals a boundary component of $N$.  It
then follows that $N'$ is a component of the exterior of $N$.   \qed

\begin{figure}[ht]
\centerline{\input{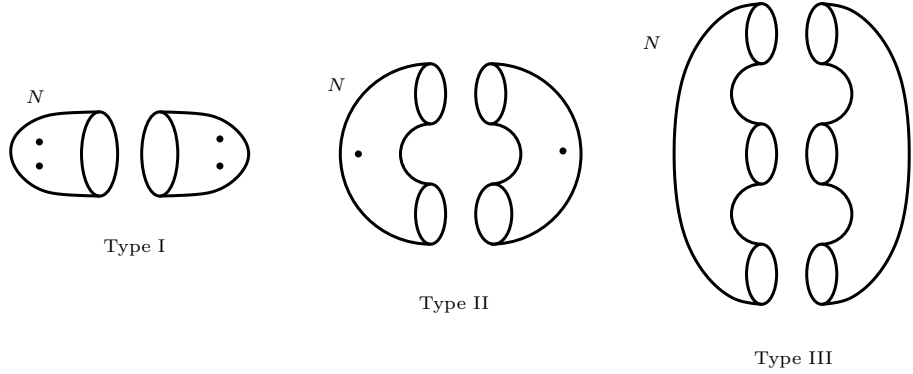}}
\caption{\label{f20} Doubled pantalons.}
\end{figure}

\bigskip\noindent
{\bf Definition:}  A {\it doubled pantalon} is the marked surface obtained by
pasting together two pantalons of the same type along their boundaries  
(see Figure \ref{f20}).
More
specifically, if $(M,P)$ is a marked surface and $N$ a subsurface of $M$, the
triple $(M,N,P)$ is called a doubled pantalon in each of the cases:

Type I: $M \cong S^2, ~~ N \cong D^2, ~~ |P| = 4, ~~ |P \cap N| = 2$;

Type II: $M \cong S^1 \times S^1 \supset S^1 \times I \cong  N, ~~ |P| = 2,  ~~
|P \cap N| = 1$;

Type III:  $M =$ closed surface of genus two, $P =$ empty set, $N \cong
\overline{M \setminus N} =$ type III pantalons.

One reason to be interested in doubled pantalons is that they provide examples
of nonisotopic subsurfaces inducing commensurable geometric subgroups, which
we will see to be an exceptional case. 

\begin{proposition} \label{p20a}  Suppose that $(M,N,P)$ is a doubled pantalon,
and let $N' = \overline{M \setminus N}$.   Then $\M(N,N \cap P)$ and 
$\M(N',N' \cap P)$ inject in $\M(M,P)$, and are commensurable subgroups.  The
subgroups are equal in the case of Types II or III.  For Type I, the
intersection $\M(N,N \cap P) \cap \M(N',N' \cap P)$ has index two in each of  
$\M(N,N \cap P)$ and $\M(N',N' \cap P)$.  In each of the three types,
the subsurfaces $N$ and $N'$ are non-isotopic in $M$ (rel $P$). 

\end{proposition}  

\pf In the case of doubled pantalons of Type II or III, the two geometric
subgroups both equal the free abelian group generated by twists along the common
boundary of $N$ and $N'$.  This has rank 2 or 3, respectively.  For Type I,
the generator of $\M(N,N \cap P)$ is a half-twist interchanging the two points
$N \cap P$, which clearly does not belong to $\M(N',N' \cap P)$, all of whose
elements fix $N \cap P$ pointwise.  Likewise the generator of $\M(N',N' \cap P)$
is not in $\M(N,N \cap P)$.  But the squares of the generators, being a Dehn
twist along the common boundary, coincide.  The question of isotopy is clear for
types $I$ and $II$, because the surfaces $N$ and $N'$ enclose different marked
points.  In type III, the two pantalons, $N$ and $N'$ are nonisotopic, too. 
For an isotopy taking $N$ to $N'$ would take the boundary to itself, but with
orientation reversed.  But a simple homological calculation shows that if
$a,b,c$ are the (oriented) boundary curves of $N$, $a$ cannot be isotopic with
$a^{-1}, b^{-1}$ or $c^{-1}.$
\qed

\bigskip
In the following, consider two connected marked  injective
subsurfaces $(N, Q), 
\break
(N', Q') \subset
(M, P)$, it being understood that $Q = N \cap P$ and $Q' = N' \cap P$.  We will
refer to the following four statements:

a) $\M(N, Q)$ and $\M(N', Q')$ are commensurable subgroups of $\M(M,P)$;

b) $\M(N, Q) = \M(N', Q')$;

c) $(N,Q)$ and $(N',Q')$ are isotopic;

d) $(N,Q)$ is isotopic with either $(N',Q')$ or $(\overline{M \setminus N'}, 
P\setminus Q')$.

\begin{theorem} \label{p21}  Suppose that 
$N$ and $N'$ are injective subsurfaces of $M$ and that 
$(N, Q) \ne (S^1 \times I, 0) \ne (N', Q')$.  Then:

i) if $(M,N,P)$ is not a doubled pantalon, 
$(a) \Leftrightarrow (b) \Leftrightarrow (c);$

ii) if $(M,N,P)$ is a doubled pantalon of Type II or III,
$(a) \Leftrightarrow (b) \Leftrightarrow (d);$

iii) if $(M,N,P)$ is a doubled pantalon of Type I,
$(b) \Leftrightarrow (c) \Rightarrow (a) \Leftrightarrow (d)$. 

\end{theorem}  

\begin{figure}[ht]
\centerline{\input{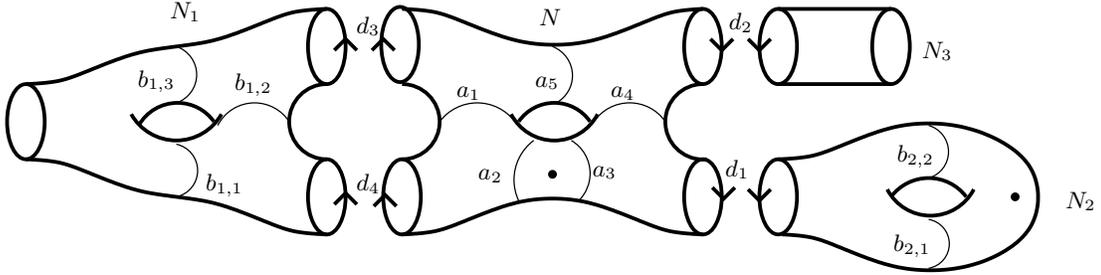}}
\caption{\label{f21} Pantalon decomposition of $N$ and its complement.}
\end{figure}

\pf  It is obvious that $(c) \Rightarrow (b) \Rightarrow (a)$ in all three
cases. Before breaking into cases, let $d_1, \dots, d_l: S^1 \to M \setminus P$ 
denote the boundary components of
$N$ and similarly $d'_1, \dots, d'_{l'}$ the components of $\partial N'$.  Let 
$\{a_i\}$ be curves in $N$ defining a pantalon decomposition of $N$
(see Figure \ref{f21}). Let $N_j$
denote the components of $\overline{M \setminus N}$.  Choose curves $b_{jk}$ giving a
pantalon decomposition of $(N_j,N_j\cap P)$
for $(N_j,N_j\cap P)\neq(S^1\times I,\emptyset)$.

We wish to show that in all cases $(a) \Rightarrow (d)$.  So we assume
commensurability of the geometric subgroups $\M(N, Q)$ and $\M(N', Q')$; the
strategy is to show that after an isotopy, the boundaries of $N$ and $N'$ can
be made to coincide.  We break the argument into two steps.

Step 1:  Each component $d'_j$ of $\partial N'$ is isotopic  (in $M$, rel $P$)
with a component $(d_i)^{\pm 1}$ of $\partial N$, and vice-versa.

Proof of Step 1: First note that we may assume in this step that
$\partial M = \emptyset$, by the trick used before: adjoin a genus one surface
to each boundary component of $M$ to obtain $\hat{M}$.  $\M(M,P)$ injects in  
$\M(\hat M,P)$, so the subgroups $\M(N,Q)$ and $\M(N',Q')$ are unchanged, and
the desired conclusion of step 1, if true in $\hat{M}$, will also hold in $M$ by
Proposition \ref{p5}.  
In particular, since $N$ is injective, none of the $N_j$ is
a cylinder with $N_j\cap P=\emptyset$, and the union of all
the $a_i, b_{jk}$ and $d_i$ gives a pantalon decomposition
of $M$.
Let
$A_i, B_{jk}$, and $D_i$ denote the Dehn twists about the respective curves
$a_i, b_{jk}, d_i$.  Let $d' = d'_j$ be any boundary component of $N'$, and
$D'$ the element of $\M(N',Q') \subset \M(M,P)$ represented by a Dehn twist
along $d'$.  We now calculate some intersection numbers:

$I(a_i, d') = 0.$ Reason:  Being a twist on a boundary curve, $D'$ is
central in $\M(N',Q')$, though perhaps not in $\M(M,P)$.  By commensurability
and Proposition \ref{p18.5} since 
$A_i \in \M(N,Q)$ there is an integer $s \ne 0$ such that 
$A_i^s \in \M(N',Q')$.  It follows that $A_i^s$ and $D'$ commute.  Then
Proposition \ref{p7} implies that $I(a_i, d') = 0$. 

$I(b_{jk}, d') = 0.$ This is shown similarly, noting that $B_{jk}$
commutes with all of $\M(N,Q)$, and by commensurability, a nonzero power of $D'$
is in $\M(N,Q)$.

$I(d_i, d') = 0$ for all $i$.  As in the previous case, a power of $D'$
belongs to $\M(N,Q)$ and so commutes with $D_i$, a twist of $\partial N$. 

Now Proposition \ref{p19} implies that $d'$ may be assumed, after an isotopy, to
be disjoint from the $a_i, b_{jk}, d_i$, i. e. $d'$ lies entirely inside
one of the pantalons in the decomposition of $M$ corresponding to those curves. 
Being essential, $d'$ is isotopic with one of the boundary curves of that
pantalon (or its inverse).  We have concluded that, up to isotopy of $M$ rel
$P$, we have one (and only one) of the following possibilities, for some
$i,j,k$:
$$d' \simeq a_i^{\pm 1}, \quad d' \simeq b_{jk}^{\pm 1},\ {\rm or} 
\quad d' \simeq d_i^{\pm 1}.$$

First assume $d' \simeq a_i^{\pm 1}$.  Since $a_i$ is 
not isotopic with a boundary component of
$N$, there exists a generic simple closed
curve $e$ in $N \setminus Q$ such that $I(d',e) = I(a_i,e) > 0$.  The
Dehn twist $E$ along $e$ is an element of $\M(N,Q)$, a subgroup
commensurable with $\M(N',Q')$, so there exists $t > 0$ such that 
$E^t \in \M(N',Q')$.  Recalling that $D'$ is central in 
$\M(N',Q')$, $D'$ and $E^t$ commute.  By Proposition \ref{p7},
$I(d',e)=0$, a contradiction.

Next assume $d' \simeq b_{jk}^{\pm 1}$, one of the pantalon curves of the
component $N_j$ of the complement of $N$.  As above, there is a simple closed
curve $e$ in $N_j \setminus (N_j \cap P)$ with  $I(d',e) = I(b_{jk},e) > 0$.
Noting that the Dehn twist $E$ commutes with all elements of $\M(N,N \cap P)$
and that a power of $D'$ belongs to $\M(N,N \cap P)$, we obtain again the
contradiction that $I(d',e)=0$.  Therefore $d' \simeq d_i^{\pm 1}$.  So Step 1
is established, noting that by symmetry, any $d_i$ is isotopic to some $d'_k$.

Step 2:  
Now, the boundary of $M$ is not assumed to be necessarily
empty.
If, as in the hypothesis, $(N,Q) \ne (S^1 \times I, 0) \ne (N', Q')$,
the curves $d_1, \dots, d_l$ are isotopically distinct, and likewise 
$d'_1, \dots, d'_{l'}$.  We conclude that
$l = l'$ and, after renumbering, each $d'_i$ is isotopic in $M$ rel $P$, with
$d_i$ or $d_i^{-1}$.  By Proposition \ref{p10} there is an isotopy of $M
\setminus P$ taking $\partial N$ to $\partial N'$.  
It follows that $(a) \Rightarrow (d)$ in all cases.

We now show that the possibility $N = \overline{M \setminus N'}$ can occur
only if $(N,Q)$ is a pantalon in which case $(M,N,P)$ is a doubled
pantalon.
Suppose that
$(N,Q)$ is not a pantalon (or $(S^1 \times I, 0), (D^2, 0), (D^2, 1)$ which are
excluded by assumption).  Then there exist curves $a, b: S^1 \to N \setminus Q$
with $I(a,b) > 0$.  Let $A, B \in \M(N, Q)$ denote the corresponding Dehn
twists.  By commensurability, $A^k \in \M(N',Q')$ for some $k > 0$.  But by
disjointness, all elements of $\M(N',Q')$ commute with all elements of
$\M(N,Q)$, so $A^k$ and $B$ commute.  Proposition \ref{p7} implies 
$I(a,b) = 0$, a contradiction.

In summary, we have shown that in case (i),  
$(c) \Rightarrow (b) \Rightarrow (a) \Rightarrow (c)$.  Cases (ii) and (iii)
also follow from the above arguments and Proposition \ref{p20a}.  \qed

\section{Normalizers and Commensurators}

The {\it commensurator} of a subgroup $H$ of a group $G$ is
$$
Com_G(H)=\{g\in G\ :\ g^{-1}Hg\ {\rm and}\ H\ {\rm are\
commensurable}\},
$$
the {\it normalizer} of $H$ in $G$ is
$$
\N_G(H)=\{g\in G\ :\ g^{-1}Hg=H\},
$$
and the {\it centralizer} of $H$ in $G$ is
$$
\ZZ_G(H)=\{g\in G\ :\ gh=hg\ {\rm for\ all}\ h\in H\}.
$$
In general we have
$$
\ZZ_G(H)\subset\N_G(H)\subset Com_G(H).\quad{\rm Also}\ H
\subset\N_G(H).
$$

For any subsurface $N$ of $M$, we define the {\it stabilizer} $Stab(N)$ as
$$Stab(N) = \{[h] \in \M(M,P) : h(N) ~{\rm is ~isotopic
~to}~ N {\rm ~in} ~M {\rm ~ rel}~ P\cup\partial M \}.$$  
Noting that $h \simeq h'$
implies $h(N) \simeq h'(N)$, we see this is a well-defined
subgroup of $\M(M,P)$.

\begin{figure}[ht]
\centerline{\input{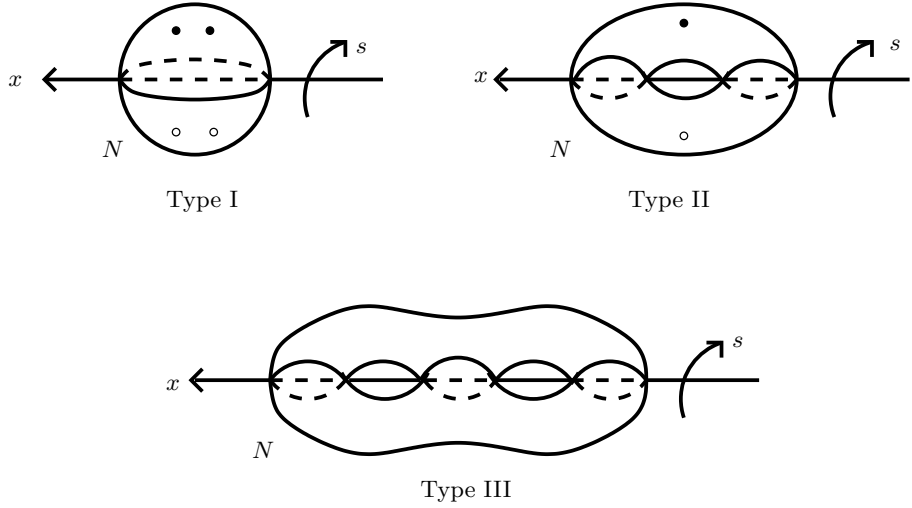}}
\caption{\label{f22} Exchange map for doubled pantalons.}
\end{figure}

As in the previous section, the doubled pantalons are exceptional cases, so
we make the following constructions.  Let $(M,N,P)$ be a doubled pantalon of
type I, II or III, embedded in the $xyz$-space symmetrically with respect to 180
degree rotation $s:M \to M$ about the $x$-axis, and so that $(N, N \cap P)$ is
the intersection of $(M, P)$ with the half-space $y \ge 0$
(see Figure \ref{f22}).  The map $s$
interchanges $N$ and $\overline{M \setminus N}$, and reverses the
orientation of $\partial N$.   We will call 
$\sigma = [s]:\M(M,P) \to \M(M,P)$, the {\it exchange} map, for each of the
three types of doubled pantalon.

\begin{theorem} \label{p22}  Let $N \subset M$ be an injective
subsurface and denote $G = \M(M, P)$ and $H = \M(N, N \cap P)$.

i) If $(M, N, P)$ is not a doubled pantalon, then
$$Com_G(H) = \N_G(H) = Stab(N).$$

ii) If $(M, N, P)$ is a doubled pantalon of
type II or III, then
$$Com_G(H) = \N_G(H) = Stab(N) \times \langle \sigma \rangle .$$

iii) If $(M, N, P)$ is a doubled pantalon of type I, then
$$Com_G(H) = Stab(N) \semidirect ~ \langle \sigma \rangle~~{\it and}~ \N_G(H) =
Stab(N).$$
Here $\langle \sigma \rangle$ denotes the cyclic subgroup of order 2 generated
by the exchange $\sigma$ of $(M,N,P)$ and $Stab(N)$ is normal in the
semidirect product. 

\end{theorem}

\pf  Let $\xi \in Com_G(H)$ and let $h \in \H(M,P)$ represent $\xi$.  Note
that 

$$\xi \M(N, P \cap N) \xi^{-1} = \M(h(N), P \cap h(N)).$$

First suppose $(M, N, P)$ is not a doubled pantalon, and if $N$ is a cylinder
assume $P \cap N$ is nonempty.  Since the groups $\M(N, P \cap N)$ and
$\M(h(N), P \cap h(N))$ are commensurable, Theorem \ref{p21} implies that
$\M(N, P \cap N) = \M(h(N), P \cap h(N))$ and that $h(N)$ is isotopic to $N$
rel $P$.  This shows that $Com_G(H) = \N_G(H) = Stab(N).$

Next suppose that $N$ is a cylinder and $P \cap N$ is empty.  Then $h(N)$ is
also a cylinder and by Proposition \ref{p20} 
$\M(N, P \cap N) = \M(h(N), P \cap h(N))$ and $h(N)$ is isotopic to $N$ rel
$P$.  Again we conclude $Com_G(H) = \N_G(H) = Stab(N).$

Now assume that $(M, N, P)$ is a doubled pantalon of type II or III.  By
Theorem \ref{p21}, $\M(N, P \cap N) = \M(h(N), P \cap h(N))$ and $h(N)$ is
isotopic to $N$ or $\overline{M \setminus N}$ rel $P$.  Therefore
$Com_G(H) = \N_G(H)$, but this group is bigger than $Stab(N)$; for example 
$\sigma$ belongs to the normalizer, but does not stabilize $N$.

We can define an epimorphism $f: Com_G(H) \to \Z /{2\Z}$ by $f([h]) = 0$ if 
$h(N)$ is isotopic to $N$ and $f([h]) = 1$ if 
$h(N)$ is isotopic to $\overline{M \setminus N}$.  The kernel of $f$ is
$Stab(N)$.  The homomorphism $f$ has a section $s: \Z /{2\Z} \to Com_G(H)$
given by $s(0) = [id], s(1) = \sigma$, where $\sigma$ is the exchange of  
$(M, N, P)$.  It follows that
$Com_G(H) = Stab(N) \semidirect ~ \langle \sigma \rangle$. 
By Propositions \ref{p14a} and \ref{p16}, $\sigma$ is in the centre of
$\M(M,P)$, so the product is, in fact, a direct product. 

Finally, consider the case that $(M, N, P)$ is a doubled pantalon of type I.
By Theorem \ref{p21}, $h(N)$ is isotopic with $N$ or 
$\overline{M \setminus N}$.  Moreover, if $\xi \in \N_G(H)$, then $h(N)$ is
isotopic with $N$, so we conclude that $\N_G(H) = Stab(N)$.  By defining the 
homomorphism $f$ and its section $s$ exactly as above, we also conclude that
$Com_G(H) = Stab(N) \semidirect ~ \langle \sigma \rangle$.  \qed

\bigskip
We proceed now to study the group $Stab(N)$ which, as seen
before, determines the commensurator and the normalizer of
$\M(N,N\cap P)$ in $\M(M,P)$.

We first state the following proposition which can be
proved in the same manner as Theorem 4.1.

\begin{figure}[ht]
\centerline{\input{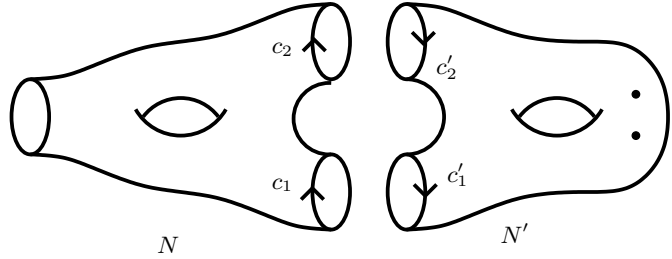}}
\caption{\label{f23} Marked surfaces whose union is $M$.}
\end{figure}

\begin{proposition} \label{p62}
Let $(N,Q)$ and $(N',Q')$ be marked surfaces. Choose $l$
boundary components $c_1,\dots,c_l$ of $N$ and $l$
boundary components $c_1',\dots,c_l'$ of $N'$, and denote
by $C_i$ and $C_i'$ the Dehn twists corresponding to $c_i$
and $c_i'$ in $\M(N,Q)$ and $\M(N',Q')$, respectively. Let
$M$ be the surface obtained by pasting together $N$ and
$N'$ along the curves $c_i$ and $c_i'^{-1}$ for all
$i=1,\dots,l$, and let $P=Q\cup Q'\subset M$ (see Figure
\ref{f23}). Assume that $|Q|\ge 2$ if $N$ is a disk and
$|Q|\ge 1$ if $N$ is an annulus. Assume the same for $N'$
and $Q'$. Consider the homomorphism 
$\phi:\M(N,Q)\times\M(N',Q')\to\M(M,P)$
induced by the inclusions of $N$ and $N'$ in $M$. Then the
kernel of $\phi$ is generated by
$$
\{(C_1,C_1'^{-1}),\dots,(C_l,C_l'^{-1})\},
$$
and is isomorphic with $\Z^l$.\qed
\end{proposition}

\bigskip
Let $(M,P)$ be a marked surface, and let $N$ be an
injective subsurface of $M$. Let $c_1,\dots,c_l$ denote
the boundary components of $N$ and $\Sigma_l$ the symmetric
group of $\{1,\dots,l\}$. There is a natural homomorphism
$\tau:Stab(N)\to\Sigma_l$ which associates to $[h]\in
Stab(N)$ the unique $\sigma\in\Sigma_l$ such that $h(c_i)$
is isotopic to $c_{\sigma(i)}$. Here we assume that the
curves $c_i$ are provided with the orientation induced by the
one of $N$. So, since $h(N)$ is isotopic to $N$, such a
$\sigma$ exists and is unique, even if $N$ is an annulus
and $N\cap P=\emptyset$. The homomorphism $\tau$ is not
surjective in general. However, one can explicitly
describe its image which depends on the topology of the
complement of $N\setminus N\cap P$ in $M\setminus P$. This
image is long and tedious to describe and we let the
reader off this description.

Let $N_1,\dots,N_r$ denote the connected components of
$\overline{M\setminus N}$. The inclusions of the $N_i$ in
$M$ and of $N$ in $M$ induce a homomorphism
$$
\phi:\M(N,N\cap P)\times\M(N_1,N_1\cap P)\times\dots\times
\M(N_r,N_r\cap P)\to\M(M,P).
$$

Assume that $N$ is not an annulus with $N\cap P=\emptyset$.
Let $C_i$ denote the Dehn twist in $\M(N,N\cap P)$
corresponding to $c_i$, $N_j$ the component of
$\overline{M\setminus N}$ having $c_i$ as boundary
component, and $C_i'$ the Dehn twist in $\M(N_j,N_j\cap
P)$ corresponding to $c_i$. Then, by Proposition 7.2, the
kernel of $\phi$ is generated by
$$
\{C_1C_1'^{-1},\dots,C_lC_l'^{-1}\}
$$
and is a copy of $\Z^l$. The image of $\phi$ is obviously
included in the kernel of $\tau$ and, by Proposition 3.10,
any element of the kernel of $\tau$ belongs in the image of
$\phi$. So, we have proved:

\begin{theorem} \label{t63}
Let $(N,N\cap P)$ be a marked injective subsurface of
$(M,P)$ different from $(S^1\times I,\emptyset)$, and let
$\tau$ and $\phi$ be the homomorphisms given above. Then we
have the exact sequence
$$
1\to\Z^l\to\M(N,N\cap P)\times\M(N_1,N_1\cap P)\times\dots 
\times\M(N_r,N_r\cap P)\stackrel{\phi}{\to}Stab(N)
\stackrel{\tau}{\to}\Sigma_l.
$$
\qed
\end{theorem}

\begin{figure}[ht]
\centerline{\input{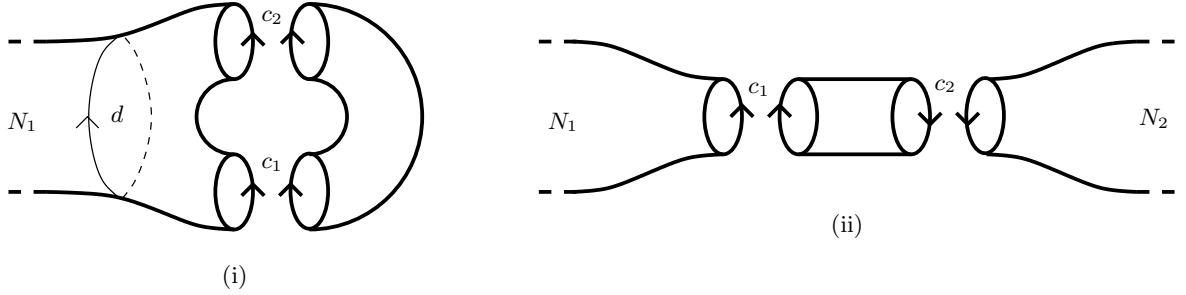}}
\caption{\label{f24} The annulus $N$ of Proposition 7.4.}
\end{figure}

\begin{proposition} \label{p64}
Assume that $N$ is an annulus, $N\cap P=\emptyset$, and
there is no boundary component of $N$ which is the boundary
of a disk intersecting $P$ in less than two points.

i) Suppose $M$ is a torus and $P$ is empty. Then
$$
Stab(N)=\M(N)\times\langle\sigma\rangle,
$$
where $\sigma$ is the element in $\M(M)$ of order two which
generates the centre.

ii) Suppose $\overline{M\setminus N}$ has a unique
connected component, $N_1$, and $(M,P)\neq(T^2,\emptyset)$
(see Figure \ref{f24}.i). Then we have an exact sequence
$$
1\to\Z\to\M(N_1,P)\to Stab(N) \stackrel{\tau}{\to}\Sigma_2
\to 1.
$$

iii) Suppose $\overline{M\setminus N}$ has two
connected components, $N_1$ and $N_2$ (see Figure
\ref{f24}.ii).Then we have an exact sequence
$$
1\to\Z\to\M(N_1,N_1\cap P)\times\M(N_2,N_2\cap P)\to
Stab(N)\stackrel{\tau}{\to}\Sigma_2.
$$
\end{proposition}

\bigskip
\pf First, suppose that $M$ is a torus and $P$ is empty. By
the structure of the mapping class group of an annulus,
the kernel of $\tau$ is $\M(N)\simeq\Z$. Moreover,
$\sigma\in Stab(N)$ and $\sigma$ permutes (up to isotopy)
the boundary components of $N$, thus $\tau$ is surjective
and the map $(1,2)\mapsto\sigma$ gives a section of
$\tau$. Since $\sigma$ is central in $\M(M)$, it
follows that
$Stab(N)=\M(N)\times\langle\sigma\rangle$. 

Assume now that $\overline{M\setminus N}$ has a unique
connected component, $N_1$, and
$(M,P)\neq(T^2,\emptyset)$. Let $c_1,c_2$ denote the
boundary components of $N$. We can cut $N_1$ along some
closed essential curve $d$ into two subsurfaces such that
one of them is a pantalon of type III having $d$, $c_1$
and $c_2$ as boundary components (see Figure \ref{f24}.i). 
Pasting this pantalon
with $N$ one obtains a genus one surface with one boundary
component, $d$. Let $\rho$ be the half-twist of this
surface of genus one along $d$ relative to $c_1$, as
defined in Section 5. Then $\rho\in Stab(N)$ and
$\tau(\rho)=(1,2)$. This shows that $\tau$ is surjective.
Let $C_i$ denote the Dehn twist along $c_i$ in
$\M(N_1,P)$, $\iota:\M(N_1,P)\to Stab(N)$ the homomorphism
induced by the inclusion of $N_1$ in $M$, and
$\theta:\Z\to\M(N_1,P)$ the homomorphism defined by
$\theta(1)=C_1C_2^{-1}$. Then the equality
$\Im\iota=\Ker\tau$ follows from Proposition 3.10, and
Theorem 4.1 implies $\Ker\theta=\{0\}$ and
$\Im\theta=\Ker\iota$.

Assume now that $\overline{M\setminus N}$ has two connected
components, $N_1$ and $N_2$. Let $c_i$ denote the common
boundary component of $N$ and $N_i$, $C_i$ the Dehn twist
along $c_i$ in $\M(N_i,N_i\cap P)$,
$$
\iota:\M(N_1,N_1\cap P)\times\M(N_2,N_2\cap P)\to Stab(N)
$$
the homomorphism induced by the inclusions of $N_1$ and
$N_2$ in $M$, and
$$
\theta:\Z\to\M(N_1,N_1\cap P)\times\M(N_2,N_2\cap P)
$$
the homomorphism defined by $\theta(1)=(C_1,C_2^{-1})$. The
equality $\Im\iota=\Ker\tau$ and the inclusion
$\Im\theta\subset\Ker\iota$ are obvious, and the
injectivity of $\theta$ follows from Corollary 3.9. So, it
remains to prove the inclusion
$\Ker\iota\subset\Im\theta$. If 
$(N_1,N_1\cap P)=(S^1\times I,\emptyset)$,
then $\M(N_1,N_1\cap P)=\M(N)$ is the infinite cyclic
subgroup generated by $C_1=C_2$, and therefore, 
since $\M(N_2,N_2\cap P)$ injects in $\M(M,P)$ by Theorem 4.1,
$\Ker\iota\subset\Im\theta$. If 
$(N_1,N_1\cap P)\neq(S^1\times I,\emptyset)\neq(N_2,N_2
\cap P)$,
then the inclusion $\Ker\iota\subset\Im\theta$ follows from
Proposition 7.2.\qed

\section{Centralizers}

\begin{theorem}\label{t71}
Let $N\subset M$ be an injective subsurface and denote
$G=\M(M,P)$ and $H=\M(N,N\cap P)$.

i) If $N$ is an annulus and $N\cap P$ is empty, then
$$
Com_G(H)=\N_G(H)=\ZZ_G(H)=Stab(N).
$$

ii) If $(M,N,P)$ is a doubled pantalon of type II or III,
then
$$
\ZZ_G(H)=\M(N,N\cap P)\times\langle\sigma\rangle,
$$
where $\langle\sigma\rangle$ is the cyclic subgroup of
order two generated by the exchange $\sigma$ of $(M,N,P)$.

iii) Suppose that $(N,N\cap P)$ is not as in (i) or (ii).
Let $N_1,\dots,N_r$ denote the connected components of
$\overline{M\setminus N}$. Then we have the exact sequence
$$
1\to\Z^l\to\ZZ\M(N,N\cap P)\times\M(N_1,N_1\cap P)\times\dots
\times\M(N_r,N_r\cap P)\stackrel{\phi}{\to}\ZZ_G(H)\to 1,
$$
where $l$ is the number of components of $\partial N$ and
$\phi$ is the homomorphism defined in Section 7. Moreover,
assuming that $|N\cap P|\ge 3$ if $N$ is a disk and $|N\cap
P|\ge 1$ if $N$ is a genus one surface with one boundary
component, the restriction
$$
\phi:\M(N_1,N_1\cap P)\times\dots\times\M(N_r,N_r\cap P)\to
\ZZ_G(H)
$$
of $\phi$ to 
$\M(N_1,N_1\cap P)\times\dots\times\M(N_r,N_r\cap P)$
is an isomorphism.
\end{theorem}

\pf Suppose that $N$ is an annulus and $N\cap P$ is empty.
The equalities $Com_G(H)=\N_G(H)=Stab(N)$ are proved in
Theorem 7.1 and the inclusion $\ZZ_G(H)\subset\N_G(H)$ is
obvious. So, it remains to prove $Stab(N)\subset\ZZ_G(H)$.
Let $c_1,c_2$ denote the boundary curves of $N$. The Dehn
twists $C_1$ and $C_2$ along $c_1$ and $c_2$, respectively,
coincide and generate $\M(N)$. Let $\xi\in Stab(N)$ and let
$h\in\H(M,P)$ represent $\xi$. Then $h\circ c_1$ is
isotopic with $c_1$ or $c_2$, thus $\xi C_1\xi^{-1}=C_1$ or
$C_2=C_1$, therefore $\xi\in\ZZ_G(H)$.

Suppose now that $(N,N\cap P)$ is not as in (i) or (ii). We
know by Theorem 7.1 that $\ZZ_G(H)\subset\N_G(H)=Stab(N)$.
Consider the exact sequence of Theorem 7.3. Let
$c_1,\dots,c_l$ denote the boundary components of $N$ and
$C_i$ the Dehn twist in $\M(M,P)$ corresponding to $c_i$.
Let $\xi\in\ZZ_G(H)$ and let $h\in\H(M,P)$ represent
$\xi$. The transformation $h$ cannot permute the $c_i$
because $\xi C_i\xi^{-1}$ is the Dehn twist along $h\circ
c_i$, the equality $\xi C_i\xi^{-1}=C_i$ implies by
Proposition 3.6 that $h\circ c_i$ is isotopic with $c_i$ or
$c_i^{-1}$, and $(N,N\cap P)\neq (S^1\times I,\emptyset)$
together with the injectivity of $N$
imply that $c_i^{-1}$ is not isotopic with some $c_j$,
$j\neq i$. So, $\ZZ_G(H)\subset\Ker\tau=\Im\phi$. Then the
exact sequence
$$
1\to\Z^l\to\ZZ\M(N,N\cap P)\times
\M(N_1,N_1\cap P)\times\dots
\times\M(N_r,N_r\cap P)\stackrel{\phi}{\to}\ZZ_G(H)\to 1
$$
is a straightforward consequence of the exact sequence of
Theorem 7.3.

Suppose now in addition that $|N\cap P|\ge 3$ if $N$ is a
disk and $|N\cap P|\ge 1$ if $N$ is a genus one surface
with one boundary component. Then, by Theorem 5.6, the
centre of $\M(N,N\cap P)$ is the free abelian group of rank $l$
generated by $\{C_1,\dots, C_l\}$. Thus it follows from
the expression of the kernel of $\phi$
given in the
proof of Theorem 7.3 that the restriction of $\phi$ to
$\M(N_1,N_1\cap P)\times\dots\times\M(N_r,N_r\cap P)$
is an isomorphism.

Assume now that $(M,N,P)$ is a doubled pantalon of type II
or III. Let $\xi\in\ZZ_G(H)$ and let $h\in\H(M,P)$
represent $\xi$. We show as above that, if $h$ stabilizes
$N$ (up to isotopy), then $h$ fixes (up to isotopy) the
boundary components of $N$, and therefore 
$h\in\M(N,N\cap P)=\M(\overline{M\setminus N},\overline{M
\setminus N}\cap P)$. So,
$$
\ZZ_G(H)\cap Stab(N)=\M(N,N\cap P).
$$
The exchange $\sigma$ is an element of the centre of
$\M(M,P)$ and thus an element of $\ZZ_G(H)$. Finally, from
the equality $\N_G(H)=Stab(N)\times\langle\sigma\rangle$ of
Theorem 7.1 follows the equality $\ZZ_G(H)=\M(N,N\cap
P)\times\langle\sigma\rangle$. \qed

\bigskip\bigskip
\parbox[t]{7cm}{
Luis PARIS\\
Laboratoire de Topologie \\
Universit\'e de Bourgogne \\
UMR 5584 du CNRS\\ 
BP 47870\\
21078 Dijon Cedex \\ 
FRANCE  \\ 
lparis@u-bourgogne.fr}
\ \ \
\parbox[t]{7cm}{
Dale ROLFSEN\\
Department of Mathematics\\
University of British Columbia\\
Vancouver\\
British Columbia V6T 1Z2\\
CANADA\\
rolfsen@math.ubc.ca}

\end{document}